\documentclass[12pt, twoside]{article}
\usepackage{amsmath,amsthm,amssymb}
\usepackage{times}
\usepackage{enumerate}
\usepackage{changes}
\usepackage{amsfonts,enumerate,graphicx,listings,epstopdf,amsthm,braket,relsize,lipsum,multirow,subcaption,chngcntr,float,wrapfig,slashbox}
\usepackage[parfill]{parskip}
\usepackage{xcolor}
\usepackage{mdframed}
\usepackage{mathtools}
\usepackage{rotating}
\usepackage{placeins}

\theoremstyle{definition}

\numberwithin{equation}{section}

\newcommand{\keywords}[1]{\noindent\emph{Keywords:}\enspace#1}

\begin{document}


\baselineskip=17pt


\title{Isogeometric Analysis of the Gray-Scott Reaction-Diffusion Equations for Pattern Formation on Evolving Surfaces and Applications to Human Gyrification}

\author{Jochen Hinz, Joost van Zwieten, Matthias M\"oller, Fred Vermolen\\
Delft Institute of Applied Mathematics\\
Delft University of Technology}

\date{}

\maketitle


\begin{abstract}
We propose a numerical scheme based on the principles of Isogeometric Analysis (IgA) for a geometrical pattern formation induced evolution of manifolds. The development is modelled by the use of the Gray-Scott equations for pattern formation
in combination with an equation for the displacement of the manifold. The method forms an alternative to the classical finite-element method. Our method is based on partitioning the initially spherical geometry into six patches, which
are mapped onto the six faces of a cube. Major advantages of the new formalism are the reconstruction of the manifold based on bicubic spline-functions used for the computation of the concentrations as well as the evolution of the mapping operator. These features give a smooth representation of the manifold which, in turn, may lead to more realistic results. The method successfully reproduces the smooth but complicated geometrical patterns found on the surface of human brains.

\end{abstract}
\newpage

\section{Introduction}
\label{chap:Introduction}
In this section, we first present the problem statement and motivation, which is followed by the introduction of the notation
and a recap of differential operators over surfaces as well as a short introduction to \textit{Isogeometric Analysis}.

\subsection{Problem statement}
The study of pattern formation is an active topic of reasearch that deals with qualitative and quantitative models for reproducing the many patterns found in the natural world, for instance in hair follicles \cite{nagorcka1985role}, leaves \cite{sachs1978patterned}, butterfly wings and mammalian coat markings \cite{murray1981pattern}. Since Turing's paper from $1952$ \cite{turing1990chemical}, a variety of reaction-diffusion (RD) systems have served as a standard model for how these patterns may arise naturally. So-called \textit{Turing Patterns} formed in the concentration of the substrates subject to these RD-type systems show a great degree of resemblance with the stripes and spots as they appear in the many labyrinth-like patterns found in nature. \\
Pattern formation, resulting from RD-models, has been verified both experimentally \cite{lee1993pattern} and computationally \cite{pearson1993complex} and the field has broadened, greatly facilitated by computational simulation. A theoretical framework for the extension of RD-models to evolving surfaces, i.e., systems in which the two-dimensional geometry is a function of time, is presented in \cite{plaza2004effect}. The evolution of the surface, in particular, may hereby be governed by the substrates on the surface. \\
Unfortunately, the underlying equations of such evolving RD-systems are often too complex to allow for an analytical treatment. This suggests the use of a numerical approach, typically one based on the principles of \textit{finite element analysis} (FEA). The purpose of the present paper is to give an example of the discretization of the Grey-Scott RD-model on an evolving geometry, initially given by a spherical shell using the principles of \textit{Isogeometric Analysis} \cite{hughes2005isogeometric}, a variant of classical FEA. To this end, we adopt the Grey-Scott RD-based model for human brain development, proposed by Lef\`evre et al. in \cite{lefevre2010reaction} and present an IgA-based numerical scheme to tackle the underlying nonlinear equations. In the following, we present a short recap of models for human brain development and a motivation for the use of IgA as a numerical technique. \\
We note that the current model provides a description of human gyrification that is too simplistic from a biophysical point of view. As such, the presented work should be considered as a feasibility study for simulating pattern formation induced geometrical development of manifolds by the use of isogeometric analysis.

\subsection{Models for Human Brain Development}
\begin{figure}[h!]
\centering
\includegraphics[width = 0.4 \linewidth]{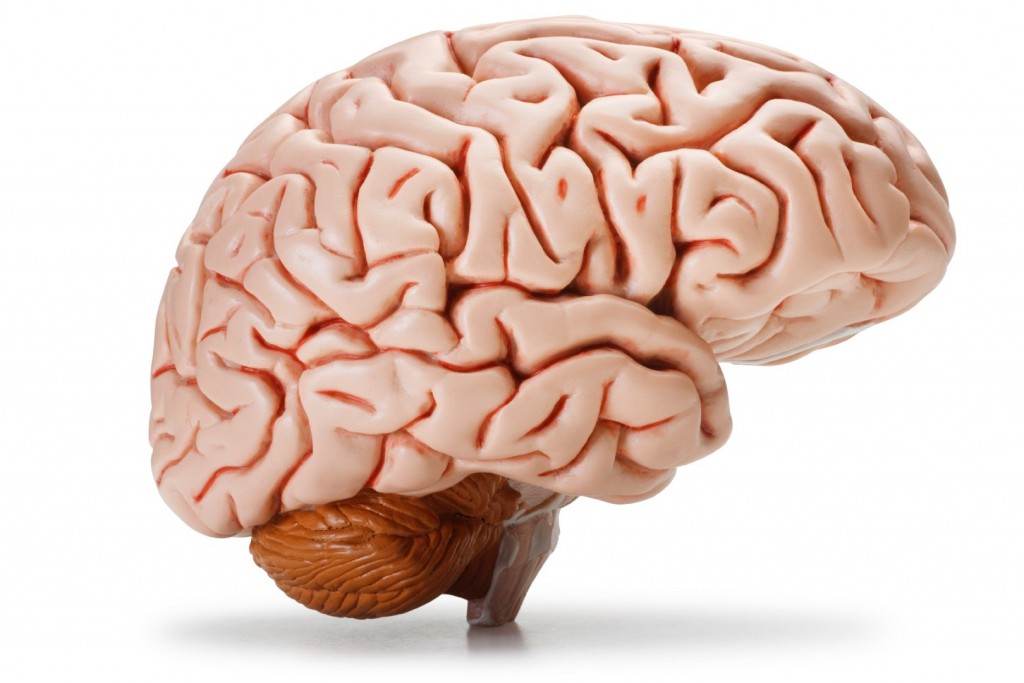}
\caption{Typical patterns formed at the surface of human brains.}
\label{fig:Human_Brain}
\end{figure}

\begin{figure}[h!]
\centering
  \begin{subfigure}[b]{0.25\linewidth}
    \centering
    \includegraphics[width=\textwidth]{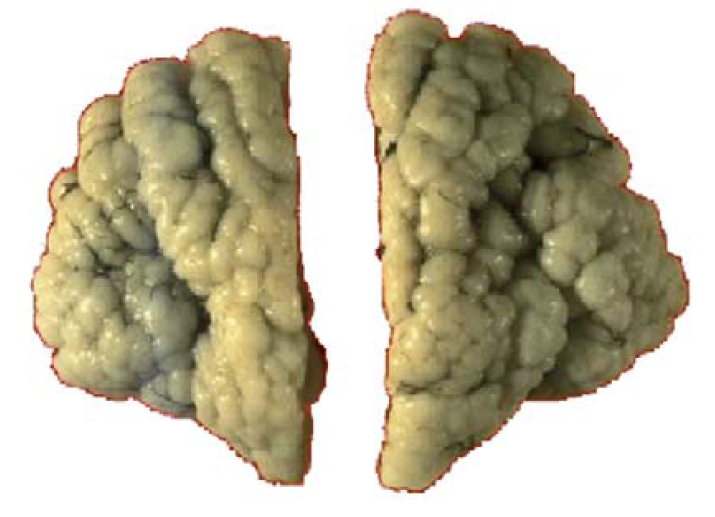}
    \caption{}
    \label{}
  \end{subfigure} \quad
  \begin{subfigure}[b]{0.25 \linewidth}
    \centering
    \includegraphics[width=\textwidth]{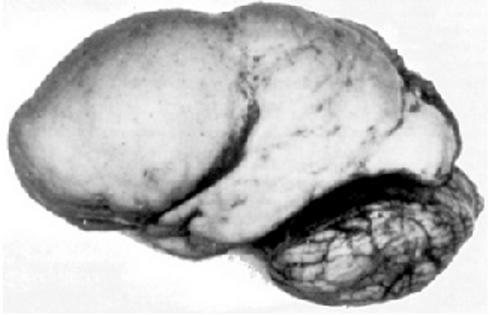}
    \caption{}
  \end{subfigure}
  \caption{Polymicrogyria (a) and Lissencephaly (b).}
\label{fig:Brain_Anomalies}
\end{figure}

Neural development has become a topic of growing interest in the past decades. On the one hand healthy adult individuals exhibit qualitatively similar neural structures, on the other hand neural development exhibits a substantial degree of randomness, which is largely confirmed by the observation that even monozygotic twins exhibit significant anatomical differences \cite{biondi1998brains}. Among other factors, this neural `fingerprint' manifests itself mainly through the patterns formed in the neural folding and buckling process occurring naturally after the twentieth week of fetal development (see Figure \ref{fig:Human_Brain}).\\
\noindent This suggests that environmental factors can have a profound influence on the course of neural development, which in turn suggests that the underlying biological process, mathematically, exhibits a high degree of sensitivity toward perturbations in the initial condition. On the other hand, a proficient model for human brain development should be capable of producing qualitatively similar outcomes for similar setups and explain neural pathologies like lissencephaly \cite{stewart1975lissencephaly} and polymicrogyria \cite{barkovich1999syndromes} (see Figure \ref{fig:Brain_Anomalies}) by quantitatively different starting conditions. \\
\noindent The derivation of proficient models for human brain development is greatly hindered by the unethicalness of experimentation on human fetuses. For this reason existing models postulate various driving forces behind pattern formation and assess their validity by comparing the results of simulations to existent brains. Furthermore quality is assessed through above criteria. Existing models for brain pattern formation can be split into two broad categories: intrinsic \cite{richman1975mechanical, raghavan1997continuum, toro2005morphogenetic, rakic2004genetic} and extrinsic \cite{van1997tension, geng2009biomechanisms, clark1945deformation}. As the name suggests, intrinsic models postulate an intrinsic process as the driving force behind neural buckling and folding, such as a chemical process, whereas extrinsic models postulate an external force, such as stresses exerted by the skull. Researchers in the field aim to identify the most-likely driving force by comparing various models in order to reduce the amount of experimentation needed to a minimum. \\
The intrinsic model proposed by Lef\`evre et al. in \cite{lefevre2010reaction} is the subject of this manuscript. It adopts a modified version of the Gray-Scott reaction-diffusion equations as a basic model for pattern formation. The concentration of one of the two chemical species considered in the model is postulated as the growth-activator, leading to deformations in the geometry that resemble typical folding patterns found in human brains. The implicit assumption of the model is that neural growth can, to a good approximation, be regarded as a process taking place only at the surface. The results from the numerical implementation presented in the article exhibit a high degree of qualitative similarities for similar setups as well as quantitative differences resulting from perturbations in the initial condition. The outcomes also show qualitative differences for different reaction rates and the authors were able to reproduce certain characteristics from various brain anomalies by changing the numerical values of the parameters. In the manuscript, we demonstrate that the various cases studied in \cite{lefevre2010reaction} can also be reproduced by the use of isogeometric analysis, in which we obtain a higher degree of smoothness in the manifold.

\subsection{Motivation}
\label{subsect:Motivation}
The model proposed by  \cite{lefevre2010reaction} results in a complex system of equations that cannot be solved analytically. The main challenge is the fact that the chemical species affect the local parametric properties of the geometry which in turn affect the local expressions of the differential operators acting on the concentrations of the chemical species, leading to a highly nonlinear system. Complexity is further increased by the existence of nonlinear reaction terms. The authors of \cite{lefevre2010reaction} present a numerical scheme that utilizes a finite-difference discretization in the temporal component, treating some terms implicitly and others explicitly, as well as a classical finite-element approach in the spatial components. The initial geometry is given by a triangularly tessellated spherical shell. \\
Growth is incorporated by extracting the magnitude of the velocity vector from the concentration of one of the chemical species at a triangle vertex, and a subsequent shift of its position in the direction of the local normal vector. On a tessellated surface, due to the non-smooth transition between adjacent triangles, the normal vector formally does not exist at the triangle vertices. The article does not explicitly state how the normal vector is computed but it is apparent that it is approximated by a weighed average of the normal vectors of the surrounding (planar) triangles. Within an implementation that is mainly focused on minimizing computational costs (in order to allow for a large amount of simulations), such an approach is reasonable. \\
Being able to construct smooth geometries would obviously constitute an improvement of above shortcoming since the normal vector is defined in every point on the geometry $\mathcal{M}$. To this end, the evident choice is to replace the classical finite-element approach with an approach based on isogeometric analysis (IgA). IgA aims to bridge the gap between classical FEA and the geometric modelling techniques from computer aided design (CAD). As such, the proposed approach, apart from some restrictions, allows for the reconstruction of the evolving surface by smooth B-spline basis functions of arbitrary polynomial order $p$ by using a spline-based mapping operator $\mathbf{s}: \Omega \rightarrow \mathcal{M}$. The geometry hereby inherits the regularity properties of the basis (locally up to $C^{p-1}$-regularity is possible) and the mapping operator becomes another unknown in the problem formulation. Numerically, we therefore treat it in the same way as the unknown concentrations on $\mathcal{M}$ and, by that, add more mathematical rigour to the scheme. \\
On the one hand, smoothness in the geometry will most-likely result in more appealing outcomes since non-smooth geometries are not realistic from a biological standpoint (we present more quantitative statements about accuracy in Section \ref{Sect:Properties_Scheme}). Since the plausibility of a model can only be assessed by the quality of the results, we hereby regard the spline-based reconstruction of the evolving surface, made possible by an IgA-based numerical approach, as particularly important. We introduce an IgA-based numerical scheme in the hopes of adding more credibility to the model. Furthermore, smoothness allows for a (non-discrete) measure of curvature, which can subsequently serve as a local refinement criterion. \\
On the other hand, modelling an evolving geometry that is initially given by a spherical shell, (not approximated by a triangular tesselation), is challenging and requires a computational domain $\Omega$ that is comprised of several quadrilaterals. With a polynomial order of $p>1$, the computational costs per degree of freedom (DOF) of the proposed IgA-scheme are expected to be higher than in the FEA-approach proposed in \cite{lefevre2010reaction} (unless $p=1$, in which case the costs are the same). Hence, this approach should be considered quality-oriented and may serve as inspiration for IgA-based numerical schemes of similar models for evolving surfaces.

\subsection{Notation}
We represent vectors and vector-valued functions utilizing bold-faced letters. The $n$-th entry of vector(-valued function) $\mathbf{r}$ is denoted by $r_n$ or $(\mathbf{r})_n$. Matrices are presented in square brackets. The $n$-th entry in the $m$-th column of matrix $[A]$ is denoted by $[A]_{n,m}$. \\
Let $\boldsymbol{\xi} = \left(\xi_1, \ldots, \xi_n \right)^T$ and $\mathbf{x} = \left(x_1, \ldots, x_m \right)^T$. We define the vector-by-vector derivative (Jacobian matrix) of $\mathbf{x}$ with respect to $\boldsymbol{\xi}$ as follows:
\begin{equation}
\label{eq:Jacobian_Introduction}
\frac{\partial \mathbf{x}}{\partial \boldsymbol{\xi}}  = \begin{bmatrix} \frac{\partial x_1}{\partial \xi_1} & \frac{\partial x_1}{\partial \xi_2} & \hdots & \frac{\partial x_1}{\partial \xi_n} \\
										    \frac{\partial x_2}{\partial \xi_1} & \frac{\partial x_2}{\partial \xi_2} & \hdots & \frac{\partial x_2}{\partial \xi_n} \\
										    \vdots 				    & \vdots 				         & \ddots & \vdots \\
										     \frac{\partial x_m}{\partial \xi_1} & \frac{\partial x_m}{\partial \xi_2} & \hdots & \frac{\partial x_m}{\partial \xi_n} \end{bmatrix}.
\end{equation}
\begin{wrapfigure}{r}{0.3\textwidth}
  \begin{center}
    \includegraphics[width = 1 \linewidth]{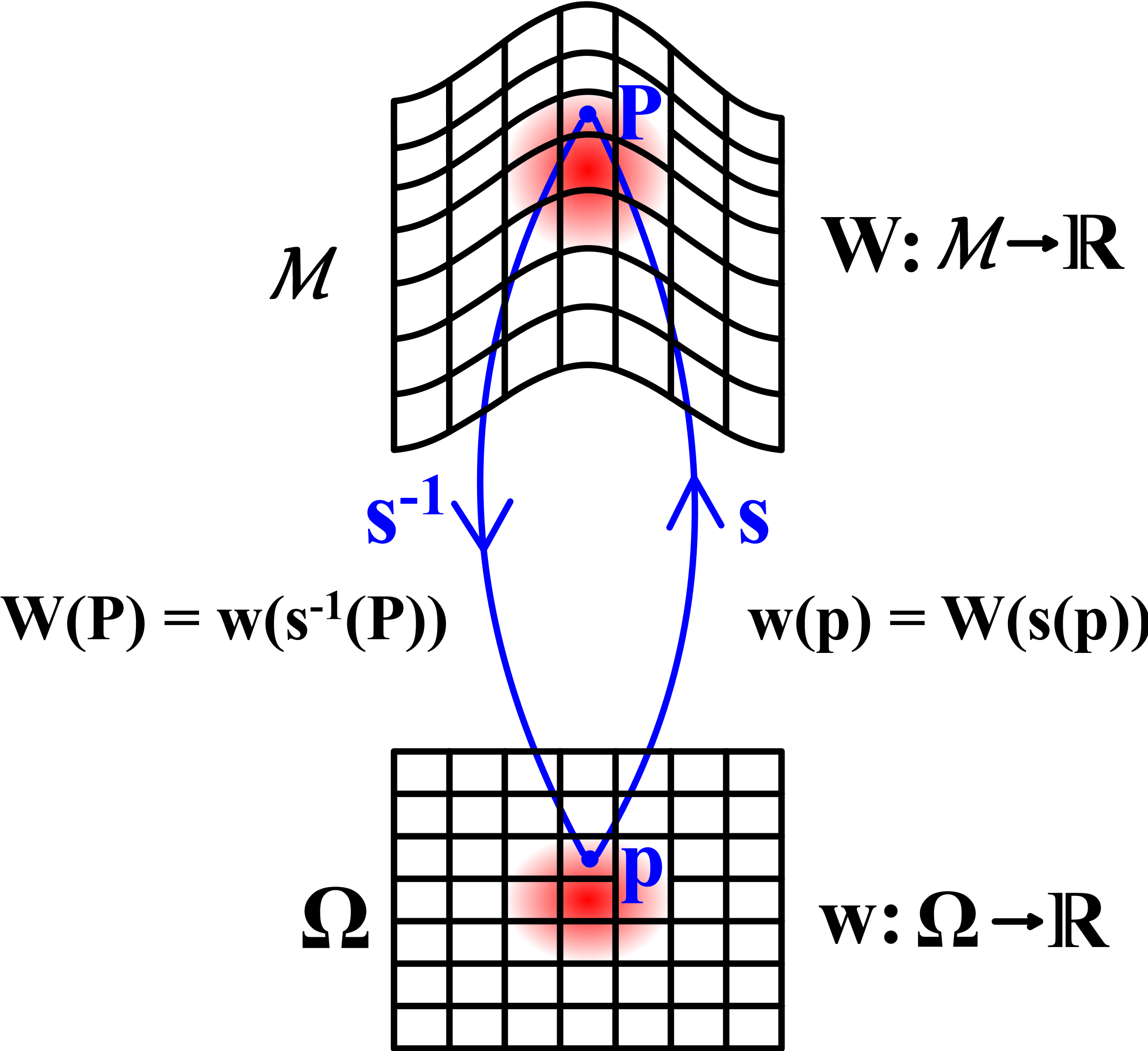}
  \end{center}
  \caption{Depiction of the relation between local and global functions.}
  \label{fig:Relation_Local_Global}
\end{wrapfigure}
In this manuscript, we will frequently work with finite-element bases. Usually the basis is denoted by 
$\Sigma = \{w_1,\ldots, w_N \}$. When referring to elements of $\operatorname{span} \Sigma$, we allow for elements with weights from $\mathbb{R}^n$ with $n \neq1$. The value of $n$ will always be clear from the context.
We make frequent use of parameterizations. By $\mathbf{s}:  \Omega \rightarrow \mathcal{M}$, with $\Omega \subset \mathbb{R}^n$ and $\mathcal{M} \subset \mathbb{R}^m$, $m \geq n$, we denote the mapping that parameterizes the geometry $\mathcal{M}$ with points from the parametric domain $\Omega$. `Local' functions living on $\Omega$ are generally represented utilizing lower-case letters, for example $w: \Omega \rightarrow \mathbb{R}$. Functions living in parameter space $\Omega$ can be made `global' by utilizing the inverse $\mathbf{s}^{-1}: \mathcal{M} \rightarrow \Omega$ of $\mathbf{s}$ and are generally represented by the corresponding upper-case letter, here: $W: \mathcal{M} \rightarrow \mathbb{R}$. Their relationship is given by
\begin{align}
\label{eq:Relation_Global_Local_Scalar_Functions}
W = w \circ \mathbf{s}^{-1},
\end{align} 
see Figure \ref{fig:Relation_Local_Global}. \\

\subsection{Differential Operators on Geometric Objects}
\label{sect:Differential_Operators_on_Manifolds}
Let
\begin{align}
[J] = \frac{\partial \mathbf{s}}{ \partial \boldsymbol{\xi} }
\end{align}
denote the Jacobian matrix of the mapping $\mathbf{s}: \Omega \rightarrow \mathcal{M}$. By
\begin{align}
[g] = [J]^T [J],
\end{align}
we denote the metric tensor associated with $\mathbf{s}$. For convenience, we introduce the short hand notation
\begin{align}
\sqrt{g} = \sqrt{\det [g]}
\end{align}
for the canonical geometric factor induced by the metric. \\
Let $\mathbf{U}: \mathcal{M} \rightarrow T_P \mathcal{M}$ and $W: \mathcal{M} \rightarrow \mathbb{R}$, where $T_P \mathcal{M}$ denotes the tangent plane of $\mathcal{M}$ spanned by the column vectors of $[J]$. Similarly to scalar functions, the local counterparts of functions $\mathbf{U}: \mathcal{M} \rightarrow T_P \mathcal{M}$ receive the corresponding lower case letter. Their relationship is hence given by:
\begin{align}
\mathbf{U} = [J] \mathbf{u}.
\end{align}
Hence, $[g]$ induces a canonical metric for vector-valued functions $U: \mathcal{M} \rightarrow T_P \mathcal{M}$ in local coordinates, since
\begin{align}
\mathbf{A} \cdot \mathbf{B} = \mathbf{a}^T [J]^T [J] \mathbf{b} = \mathbf{a}^T [g] \mathbf{b} \equiv \langle \mathbf{a}, \mathbf{b} \rangle.
\end{align}
Analogous to the standard divergence, the geometric divergence is defined as the negative adjoint of the surface gradient. This translates to local coordinates as follows:
\begin{align}
\label{eq:Definition_Geometric_Divergence}
\int_{\Omega}  w \left( \nabla_\mathcal{M} \cdot \mathbf{u} \right)  \sqrt{g} \mathrm{d} \boldsymbol{\xi} & = - \int_{\Omega} \langle \nabla_\mathcal{M} w,  \mathbf{u} \rangle \sqrt{g} \mathrm{d}  \boldsymbol{\xi},
\end{align}
for all functions $W$ that vanish on $\partial \mathcal{M}$ (compact support). The function that satisfies (\ref{eq:Definition_Geometric_Divergence}) is, in local coordinates, given by \cite[p. 18]{rosenberg1997laplacian}
\begin{align}
\label{eq:Geometric_Divergence}
\nabla_\mathcal{M} \cdot \mathbf{u} & = \frac{1}{\sqrt{g}} \sum_{i = 1}^n \frac{\partial}{\partial \xi_i} \left(\sqrt{g} u_i\right).
\end{align}
In the remainder, we shall replace $\nabla_\mathcal{M} \rightarrow \nabla$ for convenience.  \\
Similarly, the geometric gradient, in local coordinates, is given by \cite[p. 62]{ivancevic2014applied}
\begin{align}
\nabla w = [g]^{-1} \hat{\nabla} w,
\end{align}
where $\hat{\nabla}$ denotes the nabla operator in local coordinates. \\
The Laplace-Beltrami operator, being the counterpart of the ordinary Laplace-operator, consequently satisfies $\Delta w \equiv \nabla  \cdot \nabla w$. After some rearrangement, we find \cite{kreyszig1991differential}
\begin{align}
\label{eq:Laplace_Beltrami}
\Delta w = \frac{1}{\sqrt{g}} \sum_{i,j = 1}^n \frac{\partial}{\partial \xi_i} \left(\sqrt{g} g^{i,j} \frac{\partial}{\partial \xi_j} w \right),
\end{align}
where the $g^{i,j}$ are the entries of $[g]^{-1}$. \\
The above expressions can be regarded as generalizations of the commonly-encountered differential operators on manifolds. As such, most of the operations involving partial integration that are commonly applied in FEA are still applicable. In particular
\begin{align}
\label{eq:partial_integration}
\int_{\mathcal{M}} W \Delta U \mathrm{d} \mathbf{x} = \int_{\partial \mathcal{M}} W \nabla U \cdot \mathbf{N}_{\partial \mathcal{M}} \mathrm{d} l - \int_{\mathcal{M}} \nabla W \cdot \nabla U \mathrm{d} \mathbf{x},
\end{align}
where $\mathbf{N}_{\partial \mathcal{M}} \in T_P \partial \mathcal{M}$ denotes the unit outward normal vector from the tangent plane of $\partial \mathcal{M}$.

\subsection{Isogeometric Analysis}
\label{subsect:IgA}
\begin{figure}[t!]
    \centering
    \includegraphics[width = 0.75 \linewidth]{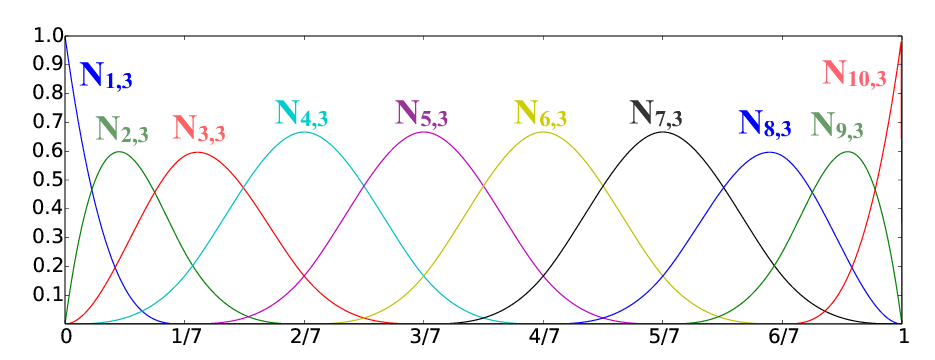}
    \caption{The univariate B-spline basis resulting from the knot vector $\Xi_3$.}
    \label{fig:open_kv_p_3}
\end{figure}

As stated in Subsection \ref{subsect:Motivation}, \textit{Isogeometric Analysis} (IgA), first introduced by Hughes et al. in \cite{hughes2005isogeometric}, is a numerical technique aimed at bridging the gap between the principles of CAD and FEA. As CAD-geometries (typically) come in the form of (B)-Spline- or NURBS-based \cite{piegl2012nurbs} mapping operators $\mathbf{s}: \Omega \rightarrow \mathcal{M}$, the \textit{isoparametric principle}, which states that the unknowns on the geometry $\mathcal{M}$ should be treated in the same way as $\mathcal{M}$ itself, lies at the heart of any IgA-based numerical simulation. To this end, the same NURBS basis used for the geometry is employed in the finite element analysis and the NURBS-based mapping operator $\mathbf{s}: \Omega \rightarrow \mathcal{M}$ is left unaltered, avoiding the need for tesselation. In the following, we present a brief recap of B-Splines, as well as spline-based surface mappings. \\
B-splines are piecewise-polynomial functions that can be constructed so as to satisfy various continuity properties at the places where the polynomial segments connect. Their properties are determined by the entries of the so-called \textit{knot vector}
\begin{align}
    \Xi = \{ \xi_1, \xi_2, \ldots, \xi_{n + p + 1} \}.
\end{align}
The knot vector is a non-decreasing sequence of parametric values $\xi_i \subset [0,1]$ that determine the boundaries of the segments on which the spline-basis is polynomial. Selecting some polynomial order $p$, the $p$-th order spline-functions $N_{i,p}$ are constructed recursively, utilizing the relation (with $\tfrac{0}{0} \equiv 0$)
\begin{align}
    N_{i,q}(\xi) = \frac{ \xi - \xi_i }{ \xi_{i+1} - \xi_i } N_{i, q-1}(\xi) + \frac{ \xi_{i+q+1} - \xi }{ \xi_{i+q+1} - \xi_{i+1} } N_{i+1, q-1}(\xi),
\end{align}
starting from
\begin{align}
    N_{i,0} = \left \{ \begin{array}{ll} 1 & \text{if } \xi_i \leq \xi \leq \xi_{i+1} \\
                                         0 & \text{otherwise} \end{array} \right. ,
\end{align}
and iterating until $q = p$. The support of basis function $N_{i,p}$ is given by the interval $\mathcal{I}_{i,p} = [\xi_i, \xi_{i+p+1}]$ and the amount of continuous derivatives across knot $\xi_j$ is given by $p - m_j$, where $m_j$ is the multiplicity of $\xi_j$ in $\mathcal{I}_{i,p}$. In practice, $\xi_1 = 0$ is repeated $p+1$ times as well as $\xi_{n+p+1}$ such that $\xi_1 = \ldots = \xi_{p+1} = 0$ and $\xi_{n+1} = \ldots = \xi_{n+p+1} = 1$. As a result, the resulting basis $\sigma = \{ N_{1,p}, \ldots, N_{n,p} \}$ forms a non-negative partition of unity on the entire parametric domain $[0,1]$, that is:
\begin{align}
    \sum_{i = 1}^n N_{i,p}(\xi) = 1,
\end{align}
with
\begin{align}
    N_{i,p}(\xi) \geq 0,
\end{align}
for all spline functions $N_{i,p}$ \cite{hughes2005isogeometric}. Figure \ref{fig:open_kv_p_3} shows the $p = 3$ B-spline basis resulting from the knot vector
\begin{align}
    \Xi = \left \{0, 0, 0, 0, \tfrac{1}{7}, \tfrac{2}{7}, \tfrac{3}{7}, \tfrac{4}{7}, \tfrac{5}{7}, \tfrac{6}{7}, 1, 1, 1, 1 \right \}.
\end{align}

\noindent The extension to bivariate spline bases is now straight-forward: given two univariate bases $\sigma_\Xi = \{N_1, \ldots, N_n\}$ and $\sigma_{\mathcal{H}} = \{M_1, \ldots M_m \}$, we build a bivariate basis $\Sigma = \{ w_{i,j} \}_{(i,j) \in \{1, \ldots, n \} \times \{ 1, \ldots, m \} }$, living on $\Omega = [0, 1]^2$, by means of a tensor-product, where
\begin{align}
    w_{i, j}(\xi, \eta) = N_i(\xi) M_j(\eta).
\end{align}
The values contained in $\Xi$ and $\mathcal{H}$ (without knot-repetitions) hereby become the boundaries of the polynomial segments, which can be regarded as the counterparts of classical elements. \\
We construct the mapping of a B-spline surface $\mathbf{s}: \Omega \rightarrow \mathcal{M}$ as follows:
\begin{align}
\label{eq:mapping_tensorial}
    \mathbf{s} = \sum_{i = 1}^n \sum_{j = 1}^m \mathbf{c}_{i,j} w_{i,j},
\end{align}
where the $\mathbf{c}_{i,j} \in \mathbb{R}^d$ are referred to as the \textit{control points}. If the dimensionaly of the control points is given by $d=3$, we are thus parameterizing a two-dimensional manifold $\mathcal{M} \subset \mathbb{R}^3$. Replacing the tensor-product index $(i, j)$ by a single global index $i$, we may therefore say that $\mathbf{s}$ is contained in $\operatorname{span} \Sigma$, with $\Sigma = \{ w_1, \ldots, w_N \}, \enskip N = nm$. The isoparametric principle now states that the same basis $\Sigma$ should be used for the unknowns on $\mathcal{M}$ parameterized by $\mathbf{s}$. As such, integrals of the form
\begin{align}
I = \int_{\mathcal{M}} W_i U \mathrm{d} \mathbf{x},
\end{align}
as they commonly arise in FEA, are computed by an appropriate \textit{pull-back} into $\Omega$. We have:
\begin{align}
I = \int_{\mathcal{M}} W_i U \mathrm{d} \mathbf{x} = \int_\Omega w_i u \sqrt{g} \mathrm{d} \boldsymbol{\xi},
\end{align}
where $\sqrt{g}$ is the geometric factor induced by $\mathbf{s}$ (see Subsection \ref{sect:Differential_Operators_on_Manifolds}),  while $(W_i, w_i)$ and $(U, u)$ are related through (\ref{eq:Relation_Global_Local_Scalar_Functions}).

\section{The Gray-Scott Reaction-Diffusion Model for Human Brain Development}
\label{chap:Model}
We use the Gray-Scott reaction-diffusion based model that was proposed by Lef\`evre et al 
\cite{lefevre2010reaction}.
Assuming that the reaction takes place on a planar surface, the Gray-Scott two-component model for human brain growth reads \cite{lefevre2010reaction}:
\begin{align}
\label{eq:Gray_Scott_Reaction_Diffusion}
\left \{ \begin{array}{ll} \partial_t U & = \quad d_1 \Delta U + F(1 - U) - UV^2 \\
\partial_t V & = \quad d_2 \Delta V  - (F+ H)V + UV^2 \end{array} \right. ,
\end{align}
subject to initial and boundary conditions. \\
Here $d_1$ and $d_2$ are the diffusion constants of $U$ and $V$ on the surface and $\Delta$ is the ordinary two-dimensional Laplace operator. Further, $F$ and $H$, respectively, denote the feeding/drainage rate and
the $\pm UV^2$-terms are reaction terms. The function $(U,V) = (1,0)$ constitutes a linearly stable steady-state solution of (\ref{eq:Gray_Scott_Reaction_Diffusion}) \cite{pearson1993complex}. \\
The (non-equilibrium) solutions of (\ref{eq:Gray_Scott_Reaction_Diffusion}) are characterized by a large number of patterns for different values of $F$ and $H$ (alternating between low and high concentration).  \\
The reaction-diffusion process need not take place on a planar geometry. The geometry $\mathcal{M}$ can, for instance, be a curved surface in $\mathbb{R}^3$ parameterized by the mapping $\mathbf{s}: \Omega \rightarrow \mathcal{M}$. As reaction-diffusion type equations commonly follow from mass-conservation considerations, the local curvature of the surface must have an effect on the system of PDEs. According to \cite{plaza2004effect}, in the presence of curvature, one has to replace $\Delta \rightarrow \Delta_\mathcal{M}$ in (\ref{eq:Gray_Scott_Reaction_Diffusion}) (see Section \ref{sect:Differential_Operators_on_Manifolds}). Thus, in the global sense, the reaction-diffusion equation reads:
\begin{align}
\label{eq:Gray_Scott_Reaction_Diffusion_Global}
\left \{ \begin{array}{ll} \partial_t U & = \quad d_1 \Delta_\mathcal{M} U + F(1 - U) - UV^2 \\
\partial_t V & = \quad d_2 \Delta_\mathcal{M} V  - (F+ H)V + UV^2 \end{array} \right., \quad \text{on } \mathcal{M},
\end{align}
subject to initial and boundary conditions. \\
The above system of PDEs is easily translated to local coordinates by replacing $U \rightarrow u$ and $V \rightarrow v$. \\
The most distinguishing feature of the model is the inclusion of surface deformation. It is assumed that species $V$ acts as a growth activator and $U$ as inhibitor. To be specific, it is suggested that the species affect the mapping operator $\mathbf{s}: \Omega \times [0, \infty) \rightarrow \mathbb{R}^3$ in the following way:
\begin{align}
\label{eq:Growth_Function}
\partial_t \mathbf{s}(\boldsymbol{\xi}, t) & = l\left(u(\boldsymbol{\xi}, t) ,v(\boldsymbol{\xi}, t) \right) \mathbf{n}(\boldsymbol{\xi}, t),
\end{align}
where $l(u,v)$ is some growth function and $\mathbf{n}$ is the unit outward normal vector to the surface. The most straightforward choice for $l(u,v)$ is
\begin{align}
\label{eq:Choice_Growth_Function}
l(u,v) & = K v,
\end{align}
for some $K > 0$, as proposed in \cite{lefevre2010reaction}. With this choice, the geometry $\mathcal{M}$ will deform at places where $V$ is nonzero, which is why $V$ acts as a \emph{growth activator}. \\
As $\mathbf{s}$, in the presence of growth, becomes a time-dependent function, the geometry will receive a time-subscript $\mathcal{M} \rightarrow \mathcal{M}_t$, where $\mathcal{M}_t$ is parameterized by $\mathbf{s}(t): \Omega \rightarrow \mathcal{M}_t$. Similarly, the surface metric receives a time-indicator $[g] \rightarrow [g_t]$ as well as the Riemannian volume form $\sqrt{g} \rightarrow \sqrt{g_t}$. \\
Obviously, the inclusion of geometry deformation has an influence on mass-conservation considerations over control surfaces. According to \cite{plaza2004effect}, in equation (\ref{eq:Gray_Scott_Reaction_Diffusion_Global}), the time-evolutions of the concentrations $U$ and $V$ have to be modified so as to contain additional terms. The terms are, in local coordinates, given by $- u \partial_t (\ln \sqrt{g_t}) $ and $- v  \partial_t (\ln \sqrt{g_t})$, respectively. With that in mind, the modified equations, in local coordinates, read:
\begin{align}
\label{eq:Gray_Scott_Reaction_Diffusion_Local_Growth}
\left \{ \begin{array}{ll} \partial_t u & =  - u \partial_t \left(\ln \sqrt{g_t} \right) +  d_1 \Delta_t u + F(1 - u) - uv^2 \\
\partial_t v & = - v \partial_t \left(\ln \sqrt{g_t} \right) + d_2 \Delta_t v  - (F+ H)v + uv^2 \end{array} \right., \quad \text{on } \Omega,
\end{align}
where $\Delta_t \equiv \Delta_{\mathcal{M}_{t}}$. \\
These additional terms ensure that the average concentrations decrease upon surface expansion and increase upon surface contraction. Since $- \partial_t \left(\ln \sqrt{g_t} \right) < 0$ whenever the surface locally expands and  $- \partial_t \left(\ln \sqrt{g_t} \right) > 0$ whenever it (locally) contracts, the modification makes intuitive sense. \\
The basic idea behind the model is that patterns formed on the surface will manifest themselves in surface deformations through the extension of the model via (\ref{eq:Growth_Function}) and (\ref{eq:Choice_Growth_Function}). With the right choice for $F$ and $H$, the patterns formed resemble the typical patterns found on the surface of human brains. 


\section{Isogeometric Implementation}
\label{chap:Isogeometric_Implementation}
In this section, we present a general numerical scheme with which the differential equation (\ref{eq:Gray_Scott_Reaction_Diffusion_Local_Growth}) from Section \ref{chap:Model} is tackled. \\
In the remainder of this section, we will assume that the Riemannian volume form $\sqrt{g_t}$ satisfies the following condition: for all $t$, there exist strictly positive constants $m_t, M_t$, with $M_t > m_t$, such that $m_t \leq \sqrt{g_t} \leq M_t$ (e.g., $[g_t]$ is \emph{`well-behaved'}). First we summarise very briefly the IgA procedure combined with the parametrisation of the spherical surface. Subsequently, we describe 
issues like time integration and refinement criterions. \\
Our initial geometry is given by a spherical shell. At first sight, it is tempting to use the classical parametrisation of a sphere on the basis of spherical coordinates. This strategy, however, leads to problems, since in particular near the poles the segments
are contracted towards the poles. This leads to singularities in the metric, which cause numerical problems such as large quadrature errors as well as more fundamental problems such as the fact that $\Sigma \ni w_i \in H^1(\mathcal{M}_t)$ can not be guaranteed due to the metric not being well-behaved. Therefore, we use an alternative approach in which the sphere is partitioned into six segments, which are mapped onto the faces of a (hollow) cube, where each face represents a section of the sphere. These segments are referred to as \emph{patches}. Each face can be mapped into the reference domain through an additional (trivial) pull-back. Given a hollow cube $\Omega$ with faces at $\pm 1$, we map points from the cube onto a spherical shell via the operator
\begin{align}
\label{eq:gaming_sphere_operator}
\tilde{\mathbf{s}}_0 = R \begin{pmatrix} 
									x \sqrt{ 1 - \frac{y^2}{2} - \frac{z^2}{2} + \frac{y^2 z^2}{3} } \\
									y \sqrt{ 1 - \frac{z^2}{2} - \frac{x^2}{2} + \frac{z^2 x^2}{3} } \\
									z \sqrt{ 1 - \frac{x^2}{2} - \frac{y^2}{2} + \frac{x^2 y^2}{3} }
																						\end{pmatrix},
\end{align}
where $R$ denotes the radius of the sphere.\\
The mapping from equation (\ref{eq:gaming_sphere_operator}) constitutes the initial parameterization and is part of the initial condition. The discretized initial condition, including expressions for the initial concentrations of $U$ and $V$ are expressed as elements from $\operatorname{span} \Sigma$ through least-squared projections, where the choice of the basis $\Sigma$ is discussed in Section \ref{sect:Implementation_Boundary_Conditions}.

\subsection{Temporal Discretization}
\label{sect:Temporal_Discretization}
We can formulate the system comprised of substrates $U$ and $V$ and geometry $\mathcal{M}_t$ parameterized by $\mathbf{s}(\boldsymbol{\xi}, t)$ as a system of equations in local coordinates
\begin{align}
\label{eq:System_of_Differential_Equations}
 \partial_t \begin{bmatrix} u \\  v \\ \mathbf{s} \end{bmatrix} & = \begin{bmatrix} - u \partial_t \left(\ln \sqrt{g_t} \right) + d_1 \Delta_{t} u + F(1 - u) - uv^2 \\ 
												    - v \partial_t \left(\ln \sqrt{g_t} \right) + d_2 \Delta_{t} v - (F + H) v + uv^2 \\ 
												    K v \mathbf{n} \end{bmatrix},
\end{align}
with
\begin{align}
\mathbf{n}(\mathbf{s}) & = \frac{1}{\left \| \frac{\partial \mathbf{s}}{\partial \xi} \times \frac{\partial \mathbf{s}}{\partial \eta} \right \|} \left( \frac{\partial \mathbf{s}}{\partial \xi} \times \frac{\partial \mathbf{s}}{\partial \eta} \right).
\end{align}
Defining
\begin{align}
\mathbf{u} & = \begin{bmatrix} u & v & \mathbf{s} \end{bmatrix}^T
\end{align}
and
\begin{align}
\mathbf{f}(\mathbf{u}) & = \begin{bmatrix} - u \partial_t \left(\ln \sqrt{g_t} \right) + d_1 \Delta_{t} u + F(1 - u) - uv^2 \\
												    - v \partial_t \left(\ln \sqrt{g_t} \right) + d_2 \Delta_{t} v - (F + H) v + uv^2 \\
												    K v \mathbf{n} \end{bmatrix},
\end{align}
the system of equations can be written as:
\begin{align}
\label{eq:Compact_System_of_Equations}
\partial_t \mathbf{u} & = \mathbf{f} (\mathbf{u}).
\end{align}
The individual components of $\mathbf{f}(\mathbf{u})$ shall be referred to as $f_u, f_v$ and $f_\mathbf{s}$. Note that $f_\mathbf{s}$ is itself a vector-valued function. \\
We integrate both sides of (\ref{eq:Compact_System_of_Equations}):
\begin{alignat}{3}
\label{eq:Integral_Form_System_of_Equations}
& \int_{t^k}^{t^{k+1}} \partial_t \mathbf{u} \mathrm{d} t & {}  = & \int_{t^k}^{t^{k+1}} \mathbf{f}(\mathbf{u}) \mathrm{d} t \nonumber \\
\implies \quad  & \mathbf{u} \left( t^{k+1} \right) - \mathbf{u} \left( t^k \right)  & = & \int_{t^k}^{t^{k+1}} \mathbf{f}(\mathbf{u}) \mathrm{d} t.
\end{alignat}
Before proceeding to the discretization, we define $\mathbf{u}^k$ as the approximation of $\mathbf{u} \left( t^k \right)$, resulting from the discretized scheme. \\
The right-hand-side integral of (\ref{eq:Integral_Form_System_of_Equations}) is approximated by a mixed implicit / explicit (IMEX) quadrature. Defining $h_k = t^{k+1} - t^k$, we utilize:
\begin{align}
\label{eq:Temporal_Discretization_Integral_RHS}
\int_{t^k}^{t^{k+1}} \mathbf{f}(u,v,\mathbf{s}) \mathrm{d} t & \simeq h_k \left[ \mathbf{g}(u^{k+1}, v^{k+1}, \mathbf{s}^{k}, \mathbf{s}^{k-1}) + \mathbf{h}(u^{k}, v^{k}, v^{k+1}, \mathbf{s}^k) \right],
\end{align}
Here,
\begin{align}
\mathbf{g}(u^{k+1}, v^{k+1}, \mathbf{s}^k, \mathbf{s}^{k-1}) & = \begin{bmatrix}  \begin{aligned} - & u^{k+1} \partial_t^h \left(\ln \sqrt{g_k} \right) +  d_1 \Delta_{k} u^{k+1} - F u^{k+1} \\
												    - & v^{k+1} \partial_t^h \left(\ln \sqrt{g_k} \right) + d_2 \Delta_{k} v^{k+1} - (F + H) v^{k+1} \end{aligned} \\
												     \mathbf{0}  \end{bmatrix}
\end{align}
and
\begin{align}
\label{eq:Right_Hand_Side_Function_Temporal}
\mathbf{h}(u^{k}, v^{k}, v^{k+1}, \mathbf{s}^k) & = \begin{bmatrix} \begin{aligned} - & u^k \left(v^k \right)^2 + F \\  & u^k \left(v^k \right)^2 \end{aligned} \\ K \mathbf{n}^k v^{k+1} \end{bmatrix},
\end{align}
with $\mathbf{n}^k = \mathbf{n} (\mathbf{s}^k)$, $\Delta_k \equiv \Delta_{t^k}$ and $\sqrt{g_k} \equiv \sqrt{g_{t^k}}$. \\
The expression $\partial_t^h \left(\ln \sqrt{g_k} \right)$ represents the time-discretization of $\partial_t \left(\ln \sqrt{g_k} \right)$:
\begin{align}
\partial_t^h \left(\ln \sqrt{g_k} \right) & = \frac{\ln \sqrt{g_k} - \ln \sqrt{g_{k-1}}}{h_{k-1}}.
\end{align}
We utilize a backward-difference scheme to avoid having to include $\sqrt{g_{k+1}}$ which is unknown at time-instance $t = t^k$. \\
We cast the discretized system into the form
\begin{align}
\label{eq:Implicit_Time-Step_Equation}
\mathbf{L} (\mathbf{u}^{k+1}) & = \mathbf{u}^k + h^{k} \mathbf{h} \left( u^k, v^k, v^{k+1} \right),
\end{align}
where
\begin{align}
\label{eq:Big_L_Operator}
\mathbf{L} \left( \mathbf{u}^{k+1} \right) & = \begin{bmatrix} \left \{ h_k \left[ \partial_t^h \left(\ln \sqrt{g_k} \right) -  d_1 \Delta_{k} + F \right ] + 1 \right \} u^{k+1} \\
												      \left\{ h_k \left[\partial_t^h \left(\ln \sqrt{g_k} \right) -  d_2 \Delta_{k} + F + H \right] + 1 \right \} v^{k+1} \\
												     \mathbf{s}^{k+1} \end{bmatrix}.
\end{align}
As before, the individual components of $\mathbf{L}, \mathbf{g}$ and $\mathbf{h}$ are referred to utilizing $u$, $v$ and $\mathbf{s}$ subscripts. Note that $\mathbf{L}$ is linear in ${\bf u}$,
 and that each component of $\mathbf{L} (\mathbf{u}^{k+1})$ only depends on the corresponding component in $\mathbf{u}^{k+1}$.

\subsection{Spatial Discretization}
\label{sect:Spatial_Discretization}
We first multiply equation (\ref{eq:Implicit_Time-Step_Equation}) by a test-function $W: \mathcal{M}_k \rightarrow \mathbb{R}$ and integrate over $\mathcal{M}_k$. For the sake of brevity, we utilize $\mathrm{d} \Omega_k \equiv \sqrt{g_k} \mathrm{d} \boldsymbol{\xi}$. Note that integrals over $\mathcal{M}_k$ are locally reformulated in terms of integrals over the hollow cube $\Omega$. The additional pull-back from the faces of $\Omega$ into the reference domain is omitted for the sake of simplicity. We have:
\begin{align}
\label{eq:Weak_Form_u}
& \int_{\mathcal{M}_k} W \mathbf{L}_u (U^{k+1}) \mathrm{d} \mathbf{x} = \int_{\mathcal{M}_k} W \left[ U^k + h_k \mathbf{h}_u \left (U^k, V^k \right) \right] \mathrm{d} \mathbf{x} \nonumber \\
\iff \quad & \begin{aligned}
 h_k d_1 \int_\Omega \langle \nabla_k w , \nabla_k u^{k+1} \rangle_{g_k}  \mathrm{d} \Omega_k  & + \int_\Omega w \left[ \left(\partial_t^h \left(\ln \sqrt{g_k} \right) + F \right) h_k + 1 \right] u^{k+1} \mathrm{d} \Omega_k \\
 & = \int_\Omega w \left[ u^k + h_k \mathbf{h}_u \left (u^k, v^k \right) \right] \mathrm{d} \Omega_k,
\end{aligned}
\end{align}
where we have made use of integration by parts, see equation (\ref{eq:partial_integration}), in conjunction with the assumption that $\mathcal{M}_k$ is a manifold without a boundary. \\
Similarly, we obtain:
\begin{align}
\label{eq:Weak_Form_v}
& \int_{\mathcal{M}_k} W \mathbf{L}_v (V^{k+1}) \mathrm{d} \mathbf{x} = \int_{\mathcal{M}_k} W \left[ V^k + h_k \mathbf{h}_v \left (U^k, V^k \right) \right] \mathrm{d} \mathbf{x} \nonumber \\
\iff \quad & \begin{aligned}
h_k d_2 \int_\Omega \langle \nabla_k w , \nabla_k v^{k+1} \rangle_{g_k} \mathrm{d} \Omega_k & + \int_\Omega w \left[ \left( \partial_t^h \left(\ln \sqrt{g_k} \right) + F + H \right) h_k + 1 \right] v^{k+1}  \mathrm{d} \Omega_k \\
 & = \int_\Omega w \left[ v^k + h_k \mathbf{h}_v \left (u^k, v^k \right) \right] \mathrm{d} \Omega_k.
\end{aligned}
\end{align}
As a next step, we introduce a local basis $\Sigma = \{w_1, \ldots, w_N \}$. In equations (\ref{eq:Weak_Form_u}) and (\ref{eq:Weak_Form_v}), we successively replace $w$ by $w_i, \enskip i=1, \ldots, N$ and we approximate $u^{k}, v^{k}$ and $\mathbf{s}^k$ by elements from $\operatorname{span} \Sigma$,
\begin{align}
\label{eq:Discretized_Substrates_Parametrization}
u^{k} & = \sum_{i = 1}^n c_i^k w_i, \quad v^{k} = \sum_{i = 1}^N d_i^k w_i \quad \text{and} \quad \mathbf{s}^k = \sum_{i = 1}^N \mathbf{e}_i^k w_i.
\end{align}
Introducing $[A]$, $[B]$, $[D]$, $\mathbf{f}^r$ and $\mathbf{w}$ with
\begin{align}
\label{eq:Matrices_Vectors_Weak_Form}
& [A]_{i,j} =  \int_\Omega w_i w_j \mathrm{d} \Omega_k, \quad [D]_{i,j} = \int_\Omega \langle \nabla_k w_i , \nabla_k w_j \rangle_{g_k} \mathrm{d} \Omega_k, \nonumber \\
& [B]_{i,j} = \int_\Omega w_i w_j \partial_t^h \left(\ln \sqrt{g_k} \right) \mathrm{d} \Omega_k, \quad (\mathbf{f}^r)_i = \int_{\Omega} w_i u^k \left(v^k\right)^2 \mathrm{d} \Omega_k \text{ and } (\mathbf{w})_i = \int_{\Omega} w_i \mathrm{d} \Omega_k,
\end{align}
we can construct a system of equations for $\mathbf{c}^{k+1} = \left(c_1^{k+1}, \ldots, c_{N}^{k+1} \right)^T$ and $\mathbf{d}^{k+1} = \left(d_1^{k+1}, \ldots, d_N^{k+1} \right)^T$. They satisfy:
\begin{align}
\label{eq:System_u}
\left \{ [A](1+ h_k F) + d_1 h_k [D] + h_k [B] \right \} \mathbf{c}^{k+1} & = h_k F \mathbf{w} - h_k \mathbf{f}^r + [A] \mathbf{c}^k, 
\end{align}
and
\begin{align}
\label{eq:System_v}
\left \{ [A](1+ h_k (F+H)) + d_2 h_k [D] + h_k [B] \right \} \mathbf{d}^{k+1} & = h_k \mathbf{f}^r + [A] \mathbf{d}^k.
\end{align}
Note that, strictly speaking, the matrices and vectors from equation (\ref{eq:Matrices_Vectors_Weak_Form}) should receive a time-superscript which has been omitted for the sake of brevity. \\
After systems (\ref{eq:System_u}) and (\ref{eq:System_v}) have been solved at time-instance $t^{k+1}$, we are in the position to update the geometry. Equations (\ref{eq:Implicit_Time-Step_Equation}) and (\ref{eq:Big_L_Operator}) suggest that we should update according to the relation
\begin{align}
\label{eq:mapping_update_relation}
\mathbf{s}^{k+1} & = \mathbf{s}^{k} + h_k K v^{k+1} \mathbf{n}^k.
\end{align}
We solve (\ref{eq:mapping_update_relation}) in the weak sense, leading to
\begin{align}
\int_\Omega w_i \mathbf{s}^{k+1}_j \mathrm{d} \Omega_k = \int_\Omega w_i \left( \mathbf{s}^{k} + h_k K v^{k+1} \mathbf{n}^k \right)_j \mathrm{d} \Omega_k \quad \forall (w_i, j) \in \Sigma \times \{1, 2, 3\},
\end{align}
finalizing the iteration.

\subsection{Essential Boundary Conditions and Choice of Basis}
\label{sect:Implementation_Boundary_Conditions}
\begin{figure}[h!]
\centering
\includegraphics[width = 0.25 \linewidth]{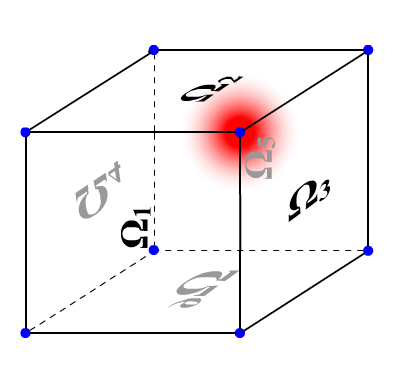}
\caption{At the cube vertices three basis functions have to be coupled (one from each patch).}
\label{fig:triple}
\end{figure}

Assuming that $\mathbf{s}^k$ parameterizes a boundaryless geometry for all $k$, there are no essential spatial boundary conditions for functions $U: \mathcal{M}_k \rightarrow \mathbb{R}$. There do, however, exist requirements for functions $u: \Omega \rightarrow \mathbb{R}$ and some initial condition
\begin{align}
\mathbf{u}(t = 0) & = \mathbf{i},
\end{align}
The projections of the components of $\mathbf{u}(t=0)$ onto the basis $\Sigma$ become the first iterates $u^0$, $v^0$ and $\mathbf{s}^0$ that jointly form the discretized initial condition. \\
The requirements we have to impose on functions $u: \Omega \rightarrow \mathbb{R}$ follow from the requirement that any $U: \mathcal{M}_k \rightarrow \mathbb{R}$ must be a function on $\mathcal{M}_k$ and hence single-valued on $\mathcal{M}_k$. In order to be compatible with the finite-element discretization, we have to furthermore require that functions $U: \mathcal{M}_K \rightarrow \mathbb{R}$ satisfy $U \in H^1(\mathcal{M}_k)$. Clearly, thanks to the well-behavedness of the metric induced by the mapping (see Section \ref{sect:Differential_Operators_on_Manifolds}), the boundedness of first-order derivatives induced by $\nabla_k$ are conditional on the boundedness of first order derivatives in local coordinates. As such, a sufficient condition for this requirement is straightforwardly translated to local functions living in $\Omega$:
\begin{align}
\label{eq:H1_requirement}
u \in H^1(\Omega) \implies U  \in H^1(\mathcal{M}_k).
\end{align}
As we shall see, this is a critical requirement especially on shared edges / vertices of the various patches under the mapping $\mathbf{s}^k$. We enforce requirement (\ref{eq:H1_requirement}) directly on the local basis functions $w_i, \enskip i = 1, \ldots, N$ such that the global counterpart $U$ of any $u \in \operatorname{span} \Sigma$ satisfies  $U \in H^1(\mathcal{M}_k)$ by default. \\
Let $\Omega_i, \enskip i = 1, \ldots, 6$ refer to the various faces of $\Omega$. The connectivities of the various patches on $\mathcal{M}_k$ follow trivially from the connectivity of the $\Omega_i$ on $\Omega$. Hence there are pairs of faces $(\Omega_i, \Omega_j)$ for which $\partial \Omega_i \cap \partial \Omega_j \neq \emptyset$. For functions $U: \mathcal{M}_k \rightarrow \mathbb{R}$, we have to impose:
\begin{align}
\label{eq:H1_requirement_patchwise}
& \left \{ \begin{array}{ll} u & \in C^0(\Omega) \\ u \vert_{\Omega_i} & \in H^1(\Omega_i), \enskip \forall i \in \{1, \ldots, 6 \}  \end{array} \right..
\end{align}
\begin{figure}[h!]
\centering
\includegraphics[width = 0.75 \linewidth]{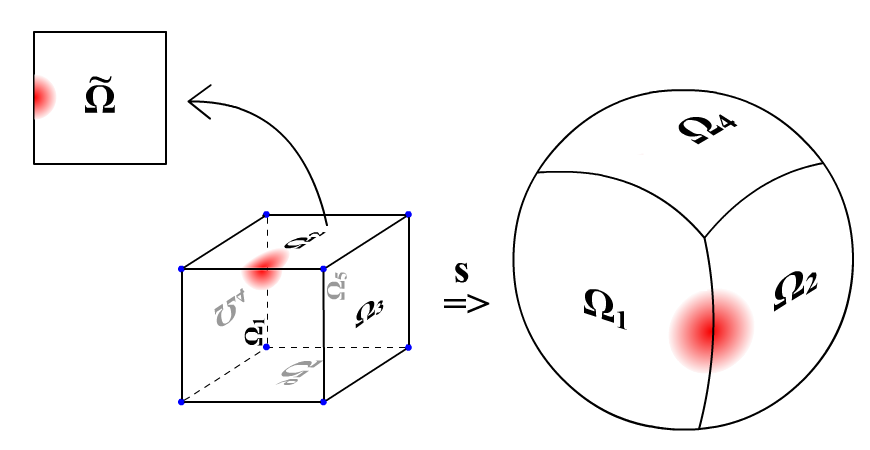}
\caption{Two functions on either side of the edge connecting $\Omega_1$ and $\Omega_2$ are combined into one function. In the absense of coupling the individual parts living on both faces would not comply with (\ref{eq:H1_requirement_patchwise}) as they are not single-valued on the connecting edge. Integrals over the cube $\Omega$ are computed by a patchwise pull back of function restrictions $u \vert_{\Omega_i}$ into the reference domain $\tilde{\Omega}$.}
\label{fig:pull_back}
\end{figure}
With condition (\ref{eq:H1_requirement_patchwise}) in mind, it becomes apparent that we can construct valid bases from the patchwise discontinuous bases $\Sigma_i, i \in \{1, \ldots, 6\}$ living on the $\Omega_i$ with appropriate degree of freedom coupling. We will assume that, upon pull-back into the reference domain, the $\Sigma_i$ constitute spline-bases with identical open and uniform knot vectors in both coordinate directions. As such, the bases are \textit{compatible} at shared edges / vertices in the sense that oppositely-faced boundary-functions can be combined into one function that is compatible with (\ref{eq:H1_requirement}) (see Figure \ref{fig:pull_back}). Besides the pairs of functions that need to be coupled on either side in the interior of shared edges, it is worth noting that the cube has eight extraordinary vertices where three functions need to be coupled (see Figure \ref{fig:triple}).

\subsection{Properties of the Numerical Scheme}
\label{Sect:Properties_Scheme}
In equation (\ref{eq:Temporal_Discretization_Integral_RHS}), we introduced the temporal discretization of the right-hand side term $\mathbf{f}(u,v,\mathbf{s})$ with
\begin{align}
\int_{t^k}^{t^{k+1}} \mathbf{f}_u \mathrm{d} t & \simeq h_k \left[ - u^{k+1} \partial_t^h \left(\ln \sqrt{g_k} \right) +  d_1 \Delta_{k} u^{k+1} - F u^{k+1} - u^k \left( v^k \right)^2 + F \right] \nonumber \\
\int_{t^k}^{t^{k+1}} \mathbf{f}_v \mathrm{d} t & \simeq h_k \left[ - v^{k+1} \partial_t^h \left(\ln \sqrt{g_k} \right) +  d_1 \Delta_{k} v^{k+1} - (F + H) v^{k+1} + u^k \left( v^k \right)^2 \right] \nonumber \\
\int_{t^k}^{t^{k+1}} \mathbf{f}_\mathbf{s} \mathrm{d} t & \simeq h_k K \mathbf{n}^k v^{k+1}.
\end{align}
The scheme can be regarded as a mixed implicit / explicit (IMEX \cite{ascher1995implicit}) scheme with the additional feature that $\mathbf{s}^{k - 1}$ is present in $\partial_t^h \left(\ln \sqrt{g_k} \right)$. The local truncation errors $\tau_k$ for the individual components of the temporal discretization are given by:
\begin{align}
\tau_k(\mathbf{f}_u) & = \mathcal{O}(h_k) + \mathcal{O}(h_{k - 1}) \nonumber \\
\tau_k(\mathbf{f}_v) & = \mathcal{O}(h_k) + \mathcal{O}(h_{k - 1}) \nonumber \\
\tau_k(\mathbf{f}_\mathbf{s})_i & = \mathcal{O}(h_k).
\end{align}
The global error is of the same order \cite{suli2003introduction}. \\
In the following, we investigate the properties of the scheme introduced in Section \ref{sect:Spatial_Discretization},  in particular symmetry and coercivity.  The weak form leads to a left hand side operator of the form:
\begin{align}
\mathbf{A}^k(W, \mathbf{U}^{k+1}) = \int_{\mathcal{M}_k} W \mathbf{L}(\mathbf{U}^{k+1}) \mathrm{d} \mathbf{x}.
\end{align}
As each component of $\mathbf{L}(\mathbf{U}^{k+1})$ only depends on the corresponding component in $\mathbf{U}^{k+1}$, the analysis can be carried out component-wisely. In equation (\ref{eq:Weak_Form_u}), we derived:
\begin{align}
\label{eq:A_u}
\mathbf{A}^k_{u}(U, V) =  h_k d_1 \int_\Omega \langle \nabla_k u , \nabla_k v \rangle_{g_k}  \mathrm{d} \Omega_k  & + \int_\Omega u \left[ \left(\partial_t^h \left(\ln \sqrt{g_k} \right) + F \right) h_k + 1 \right] v \mathrm{d} \Omega_k.
\end{align}
Clearly, $\mathbf{A}^k_u(U, V)$ is symmetric. Furthermore, we have:
\begin{align}
\mathbf{A}^k_{u}(U, U) & = h_k d_1 \int_\Omega \langle \nabla_k u , \nabla_k u \rangle_{g_k}  \mathrm{d} \Omega_k  + \int_\Omega u \left[ \left(\partial_t^h \left(\ln \sqrt{g_k} \right) + F \right) h_k + 1 \right] u \mathrm{d} \Omega_k \nonumber \\
& \geq \int_{\Omega} \left[ \left(\partial_t^h \left(\ln \sqrt{g_k} \right) + F  \right) h_k + 1 \right] u^2 \mathrm{d} \Omega_k.
\end{align}
As such, a sufficient condition for coercivity is:
\begin{align}
\label{eq:condition_coercivity}
\left(\partial_t^h \left(\ln \sqrt{g_k} \right) + F \right) h_k + 1 > 0.
\end{align}
It is easy to verify that, for a given growth rate $K$, we have:
\begin{align}
\partial_t^h \left(\ln \sqrt{g_k} \right) = \mathcal{O}(K v^{k - 1}).
\end{align}
As the growth rate $K$ typically satisfies $K \ll F$, in conjunction with the fact that the velocity vector points in the direction of the normal vector corresponding to $\mathcal{M}_k$ (leading to surface expansion, i.e., the magnitude of $\sqrt{g_t}$ is expected to increase over time), it is reasonable to assume that (\ref{eq:condition_coercivity}) is not violated. As such, $\mathbf{A}^k_u$ is coercive with coercivity constant 
\begin{align}
C^k_u = \inf_{\Omega} \left[ \left(\partial_t^h \left(\ln \sqrt{g_k} \right) + F \right) h_k + 1 \right] \sqrt{g_k}.
\end{align}
Similarly, a sufficient condition for coercivity of $\mathbf{A}^k_v$ is:
\begin{align}
\left(\partial_t^h \left(\ln \sqrt{g_k} \right) + F + H \right) h_k + 1 > 0,
\end{align}
which is likely not violated either. As such, we may conclude that $\mathbf{A}^k$ is symmetric and coercive and the fully discretized scheme symmetric positive definite and hence invertible. This suggests the use of a CG-type algorithm for the inversion of the linear systems resulting from relations (\ref{eq:System_u}) and (\ref{eq:System_v}). \\
In the following, we analyse the impact of the higher-order smoothness, made possible by the use of B-Spline basis functions, on the reconstruction of the surface that is initially given by the projection of (\ref{eq:gaming_sphere_operator}) onto $\Sigma$. For the analysis, we use a basis $\Sigma$ resulting from employing the uniform open knot vector
\begin{align}
\Xi_{n, p} = \left \{ \underbrace{0, \ldots, 0}_{p + 1 \text{ times}}, \underbrace{ \frac{1}{n - p}, \ldots, \frac{n - p -1}{n - p} }_{n - p - 1 \text{ terms}}, \underbrace{1, \ldots, 1}_{p + 1 \text{ times}}  \right \}
\end{align}
patchwise in both directions.  For fixed $n$, the use of above knot vector leads to a basis $\Sigma$ of equal cardinality, regardless of $p$.
\begin{figure}[H]
\centering
  \begin{subfigure}[b]{0.25\textwidth}
    \includegraphics[width=\textwidth]{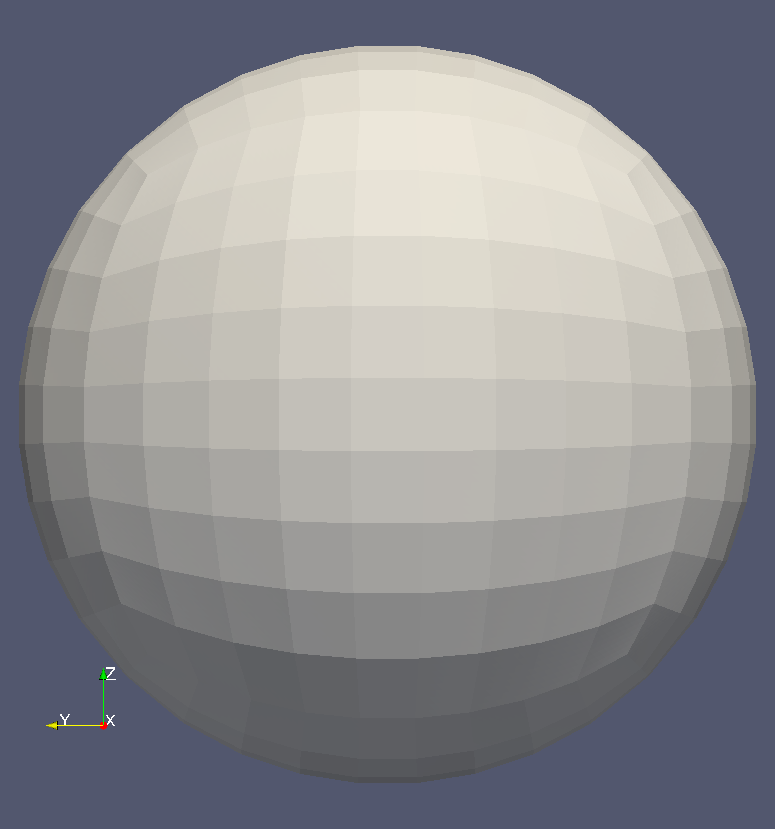}
    \caption{}
  \end{subfigure}
  \quad
  \begin{subfigure}[b]{0.25\textwidth}
    \includegraphics[width=\textwidth]{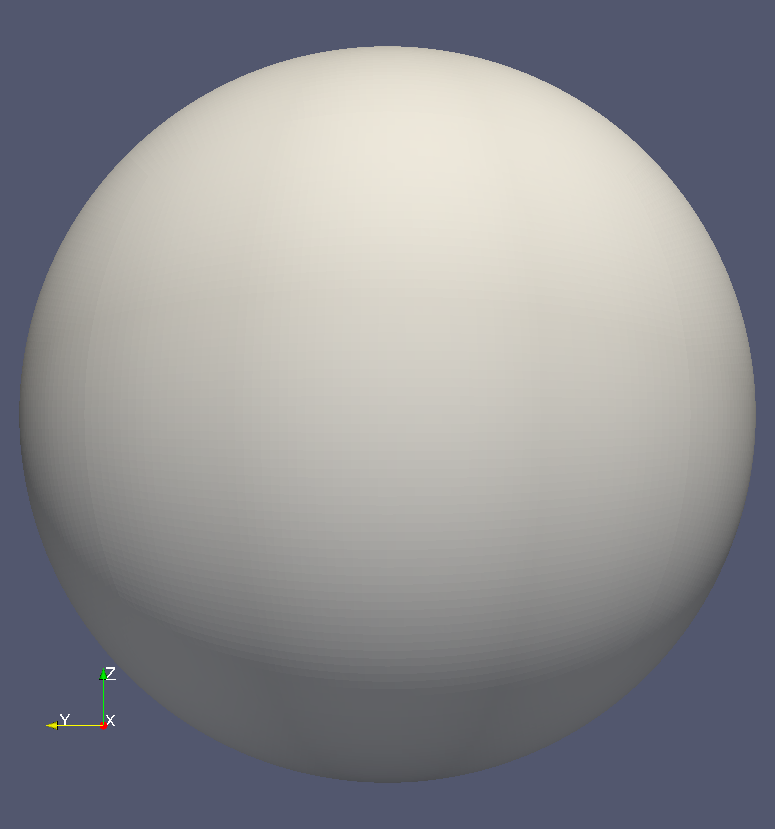}
    \caption{}
  \end{subfigure}
  \quad
  \begin{subfigure}[b]{0.25\textwidth}
    \includegraphics[width=\textwidth]{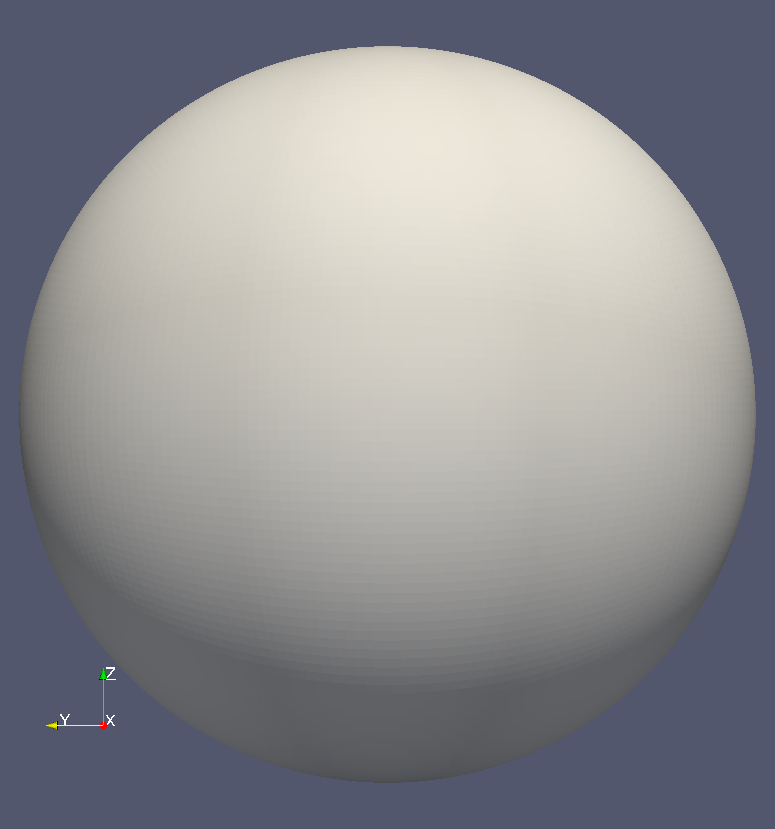}
    \caption{}
  \end{subfigure}
  \caption{Approximation of $\tilde{\mathbf{s}}_0$ using spline bases resulting from $\Xi_{10, 1}$ (a), $\Xi_{10, 2}$ (b) and $\Xi_{10, 3}$ (c). For $p=1$, we clearly see the planar elements resulting from the tesselation while the approximation becomes virtually indistinguishable from a spherical shell parameterized by $\tilde{\mathbf{s}}_0$ for $p \geq 2$.}
\label{fig:approx_sphere}
\end{figure}
\noindent Let $E_{n, p} = \left \| \tilde{\mathbf{s}}_0 - \mathbf{s}_0^{n, p} \right \|_{L^2(\Omega)}$ denote the $L^2(\Omega)$ mismatch between $\tilde{\mathbf{s}}_0$ (see (\ref{eq:gaming_sphere_operator})) and its approximation $\mathbf{s}_0^{n, p}$ resulting from $\Xi_{n, p}$.
\begin{table}[h!]
\centering
\begin{tabular}{|l||*{3}{c|}}\hline
\backslashbox{$n$}{$p$}
&\makebox[3em]{1}&\makebox[3em]{2}&\makebox[3em]{3} \\\hline\hline
$5$ & $2.0 \times 10^{0}$ & $3.6 \times 10^{-1}$ & $2.0 \times 10^{-1}$ \\\hline
$10$ & $3.8 \times 10^{-1}$ & $2.3 \times 10^{-2}$ & $4.0 \times 10^{-3}$ \\\hline
$20$ & $8.0 \times 10^{-2}$ & $2.1 \times 10^{-3}$ & $1.6 \times 10^{-4}$ \\\hline
\end{tabular}
\caption{Value of $E_{n,p}$ for various $n, p$.}
\label{tab:truncation_error_sphere}
\end{table}
Table \ref{tab:truncation_error_sphere} clearly demonstrates the superiority of the higher-order spline functions in reconstructing the spherical shell. The table shows that, for equal $n$, projections with bases of order $p$ are consistently one order of magnitude more accurate than those of order $p-1$. To achieve a similar accuracy as for $(n, p) = (10, 1)$ (amounting to $1464$ DOFs), for instance, a basis with $(n, p) = (5, 2)$ is sufficient (amounting to only $294$ DOFs). The improved accuracy is also clearly visible in Figure \ref{fig:approx_sphere}. \\
For a more in-depth analysis of the approximation properties of spline basis functions, we refer to \cite[Chapter~4]{buffa2016isogeometric}.
\subsection{Time-Step Selection}
\label{sect:Time_Step_Selection}
So far we have not discussed the choice of the time-step $h_k$. Since the system-matrices have to be rebuilt after each iteration, the time-step selection should be based on the following principles:
\begin{itemize}
\item the time-step selection should be as large as possible with respect to numerical stability;
\item the time-step selection should be small enough with respect to the characteristic time-scale of the equation to warrant numerical accuracy;
\item the time-step selection should be a cheap operation.
\end{itemize}
With the above principles in mind, time-step selection based on the proportional-integral-derivative (PID) controller proposed by Valli et al. \cite{valli2002control} constitutes a proficient choice. Here, we use the relative change of the $L_2(\mathcal{M}_k)$-norm of the components of $\mathbf{u}^k$ as control parameters at time instance $t = t^k$. This is a cheap operation since, for instance,
\begin{align}
\left \| u^k - u^{k-1} \right \|^2_{L_2(\mathcal{M}_k)} = \int_{\mathcal{M}_k}( U^k - U^{k-1} )^2 \mathrm{d} \mathbf{x} = (\mathbf{c}^k - \mathbf{c}^{k-1})^T [A] (\mathbf{c}^k - \mathbf{c}^{k-1})
\end{align}
and the mass matrix $[A]$ has to be assembled anyway. \\
The new time-step is then selected based on the following scheme:
\begin{align}
\label{eq:PID_scheme}
h_{k+1} = \left( \frac{e_{k-1}}{e_k} \right)^{k_P} 	\left( \frac{1}{e_k} \right)^{k_I} \left( \frac{e_{k-1}^2}{e_k e_{k-2}} \right)^{k_D} h_k,
\end{align}
with
\begin{align}
e_k = \max( e_u, e_v ), \quad \text{with} \quad e_u = \frac{ \left\| u^{k+1} - u^k \right \|_{L_2(\mathcal{M}_k)} }{ \tau \left \| u^{k+1} \right \|_{L_2(\mathcal{M}_k)} } \quad \text{and} \quad e_v = \frac{ \left\| v^{k+1} - v^k \right \|_{L_2(\mathcal{M}_k)} }{ \tau \left \| v^{k+1} \right \|_{L_2(\mathcal{M}_k)} }.
\end{align}
The PID-controller from equation (\ref{eq:PID_scheme})  is designed to select a time-step $h_k$ such that the relative change in the concentrations is close to, but does not exceed $\tau > 0$. As default parameters, we use $(k_P, k_I, k_D, \tau) = (0.075, 0.175, 0.01, 0.01)$.

\subsection{Refinement Strategies}
\label{sect:refinement_strategies}
In Section \ref{sect:Implementation_Boundary_Conditions}, we discussed means to construct a spline basis $\Sigma$ that is compatible with the numerical scheme. We can build a hierarchy $\{ \Sigma_1, \Sigma_2, \ldots \}$ of spline-bases with increasing cardinality by repeating the steps taken to build $\Sigma_k$ with a knot vector that results from uniformly refining the previous knot vector. It can be shown \cite[Chapter~2]{hughes2005isogeometric} that $\operatorname{span} \Sigma_1 \subset \operatorname{span} \Sigma_2 \subset \ldots $. As such, there exists a canonical sparse prolongation matrix $[\mathcal{T}_{k \rightarrow k+1}]$ that prolongs the weights of any function in $\operatorname{span} \Sigma_k$ to $\Sigma_{k+1}$. \\
Local refinement is accomplished by identifying basis functions $\Sigma_k \ni w_i^k \in \Sigma$ for refinement and replacing them by several basis functions from $\Sigma_{k+1}$.  Our refinement strategy is losely based on the principles of hierarchical B-spline refinement, see \cite{vuong2012adaptive}. Whenever some $w_i^k \in \Sigma$ has been marked for refinement, we replace it by $ \{ w_j^{k+1} \in \Sigma_{k+1} \enskip \vert \enskip [\mathcal{T}_{k \rightarrow k+1}]_{ji} \neq 0 \}$ in $\Sigma$. Basis refinement is accompanied by element refinement in order to ensure that the $w_i \in \Sigma$ remain piecewise polynomial with respect on elements of the grid. \\
We base our local refinement criterion on the properties of the mapping operator $\mathbf{s}^k$. Let
\begin{align}
\mu_k^i = \frac{ \int \limits_{\Omega} w_i \sqrt{g_k} \mathrm{d} \boldsymbol{\xi} }{ \int \limits_{\Omega} w_i \mathrm{d} \boldsymbol{\xi} }
\end{align}
and
\begin{align}
\kappa_k^i = \frac{ \int \limits_{\Omega} (\kappa_1^2 + \kappa_2^2) w_i  \sqrt{g_k} \mathrm{d} \boldsymbol{\xi} }{ \int \limits_{\Omega} w_i \sqrt{g_k} \mathrm{d} \boldsymbol{\xi} },
\end{align}
where $\kappa_1$ and $\kappa_2$ denote the two principal curvatures, defined as the two eigenvalues of 
\begin{align}
[S] = [g]^{-1} [L],
\end{align}
with
\begin{align}
[L] = \begin{bmatrix} \frac{\partial^2 \mathbf{s}}{\partial \xi^2} \cdot \mathbf{n} & \frac{\partial^2 \mathbf{s}}{\partial \xi \partial \eta} \cdot \mathbf{n} \\ \frac{\partial^2 \mathbf{s}}{\partial \xi \partial \eta} \cdot \mathbf{n} &  \frac{\partial^2 \mathbf{s}}{\partial \eta^2} \cdot \mathbf{n} \end{bmatrix}.
\end{align}
Furthermore, we define
\begin{align}
\mu_{\text{cell}} = \frac{1}{N} \sum_{i = 1}^N \mu_0^i \quad \text{and} \quad \mu_{\text{curve}} = \frac{1}{N} \sum_{i = 1}^N \kappa_0^i
\end{align}
as a measure for the initial average weighed cell size and curvature, repsectively. \\
Basis function $w_i \in \Sigma$ is marked for refinement if either of the following situations arises:
\begin{itemize}
\item $\mu_k^i > k_{\text{cell}} \mu_{\text{cell}}$;
\item $\kappa_k^i > k_{\text{curve}} \mu_{\text{curve}}$,
\end{itemize}
where $k_\text{cell} > 1$ and $k_\text{curve} > 1$ are two positive constants. \\
Refinement based on a (non-discrete) measure of curvature is made possible by the higher-order nature of the employed B-spline basis leading to a patchwise smooth parameterisation of the geometry. Since a non-discrete measure of curvature is ill-defined for basis functions that are nonvanishing on the patch interfaces (due to the local $C^0$-continuity), they are disregarded in curvature-based refinement.

\section{Results and Discussion}
\label{chap:Results}
\begin{figure}[h!]
\centering
  \begin{subfigure}[b]{0.45\textwidth}
    \centering
    \includegraphics[height=5.5cm]{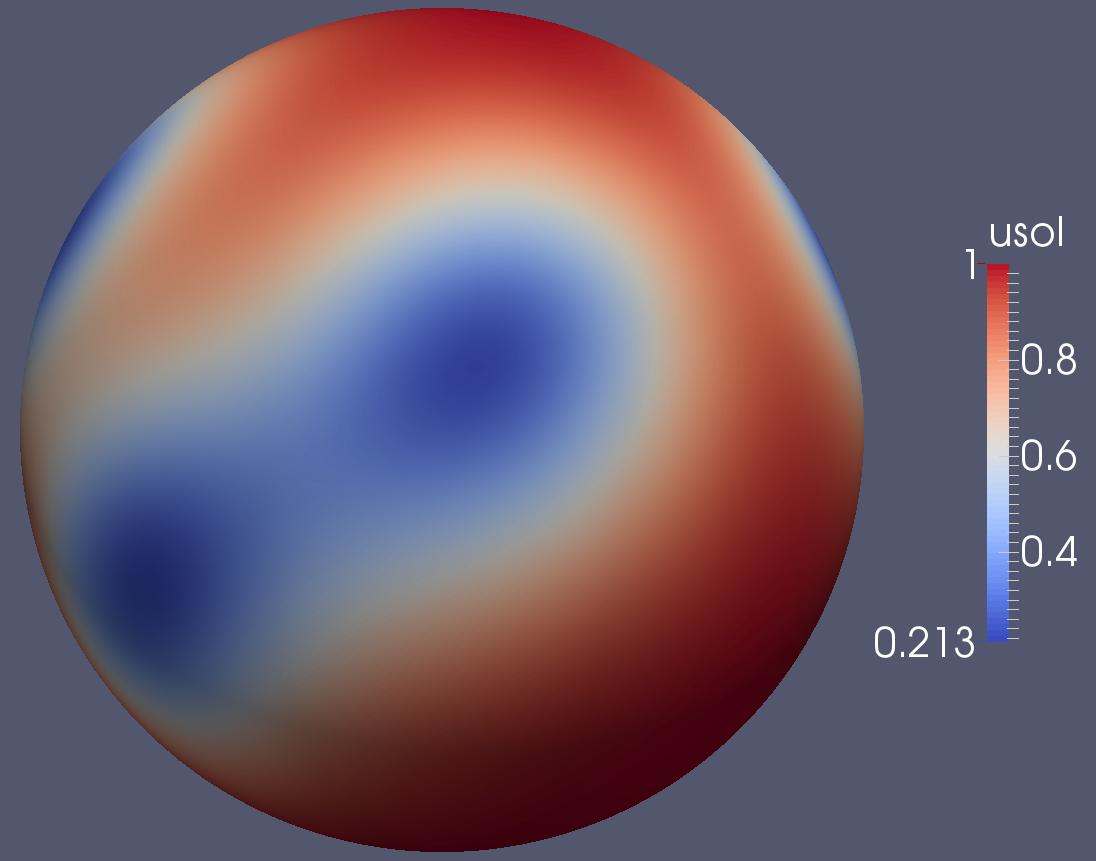}
    \caption{}
  \end{subfigure}
  \quad
  \begin{subfigure}[b]{0.45\textwidth}
    \centering
    \includegraphics[height=5.5cm]{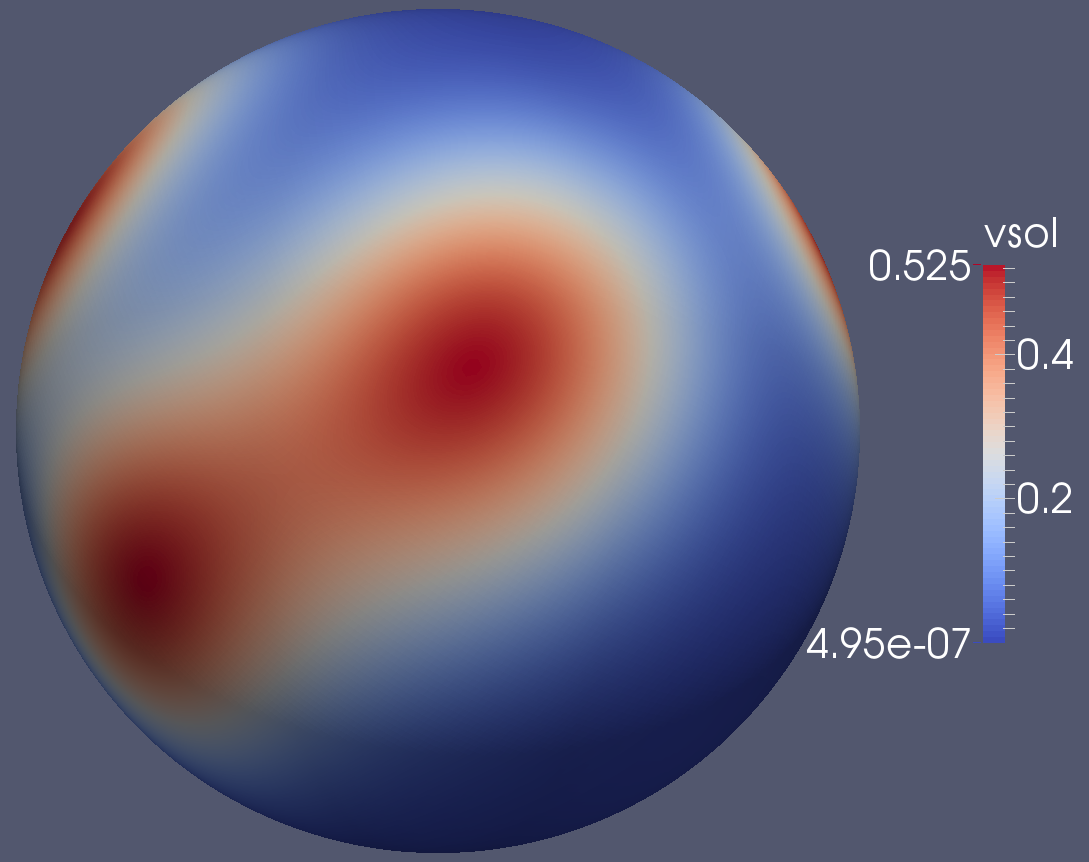}
    \caption{}
  \end{subfigure}
  \caption{Initial condition of $U$ (a) and $V$ (b) plotted on the initial geometry.}
\label{fig:Sphere_0_004}
\end{figure}

In this chapter we will present the results of an implementation with a spherical shell (see (\ref{eq:gaming_sphere_operator})) with radius $R=40$ as initial geometry. The implementation follows the principles from Section \ref{chap:Isogeometric_Implementation}, with an IgA basis based on the principles of Section \ref{sect:Implementation_Boundary_Conditions}. \\
The matrix and vector assemblies are carried out with Gauss-schemes of order six and the linear systems are solved with an iterative CG-solver. We employ a third-order B-spline basis resulting from a uniform open knot vector with $24$ internal knots in both coordinate-directions. Both refinement strategies discussed in Section \ref{sect:refinement_strategies} have been used and we employed the PID-controller from Section \ref{sect:Time_Step_Selection} for the time step selection. The implementation has been realized in the Python-package \emph{Nutils} \cite{Nutils}. \\
As in \cite{lefevre2010reaction}, the algorithm was manually terminated in order to avoid geometric self-intersections at later stages. The initial concentrations satisfy
\begin{align}
U(t = 0) & = 1 - 0.75 \times I \nonumber \\
V(t = 0) & = 0.5 \times I,
\end{align}
where the function $I$  is given by the sum of four Gaussians
\begin{align}
\label{eq:G}
G_{\mathbf{x}_0}(\mathbf{x}) & = \operatorname{exp}\left({-  \left(\frac{x_1 - (\mathbf{x}_0)_1}{20} \right)^2 -  \left(\frac{x_2 - (\mathbf{x}_0)_2}{15} \right)^2 -  \left(\frac{x_3 - (\mathbf{x}_0)_3}{15} \right)^2} \right).
\end{align} 
They are centered at
\begin{align}
\mathbf{x}_0(\xi_i,\eta_i) & = R \begin{bmatrix} \sin(\xi_i) \cos(\eta_i) \\ \sin(\xi_i)\sin(\eta_i) \\ \cos(\xi_i) \end{bmatrix},
\end{align}
with $(\xi_1,\eta_1) = (0,0)$, $(\xi_2,\eta_2) = (0.3 \pi,0.4 \pi)$, $(\xi_3,\eta_3) = (0.4 \pi,0.7 \pi)$ and $(\xi_4,\eta_4) = (0.65 \pi, \pi)$.
The initial concentrations are depicted in Figure \ref{fig:Sphere_0_004}.

\begin{figure}[h!]
\centering
  \begin{subfigure}[b]{0.45\textwidth}
    \centering
    \includegraphics[height=5.5cm]{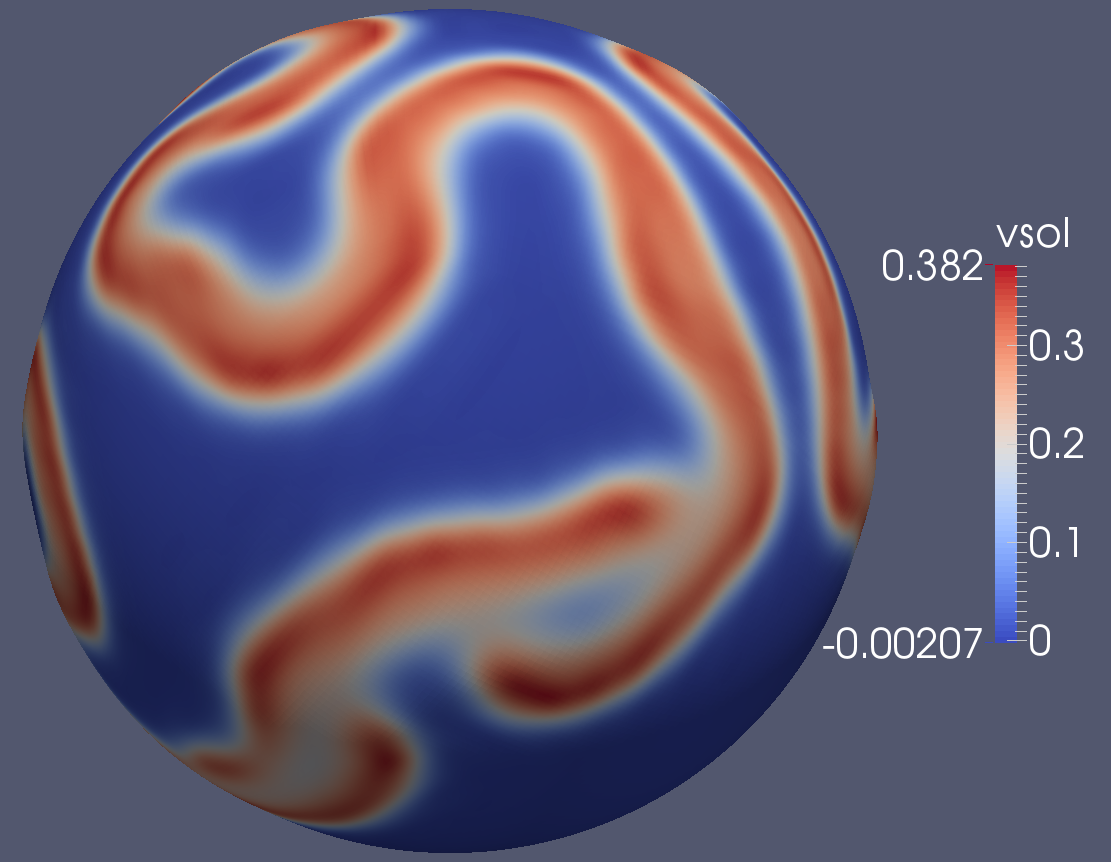}
    \caption{}
  \end{subfigure}
  \quad
  \begin{subfigure}[b]{0.45\textwidth}
    \centering
    \includegraphics[height=5.5cm]{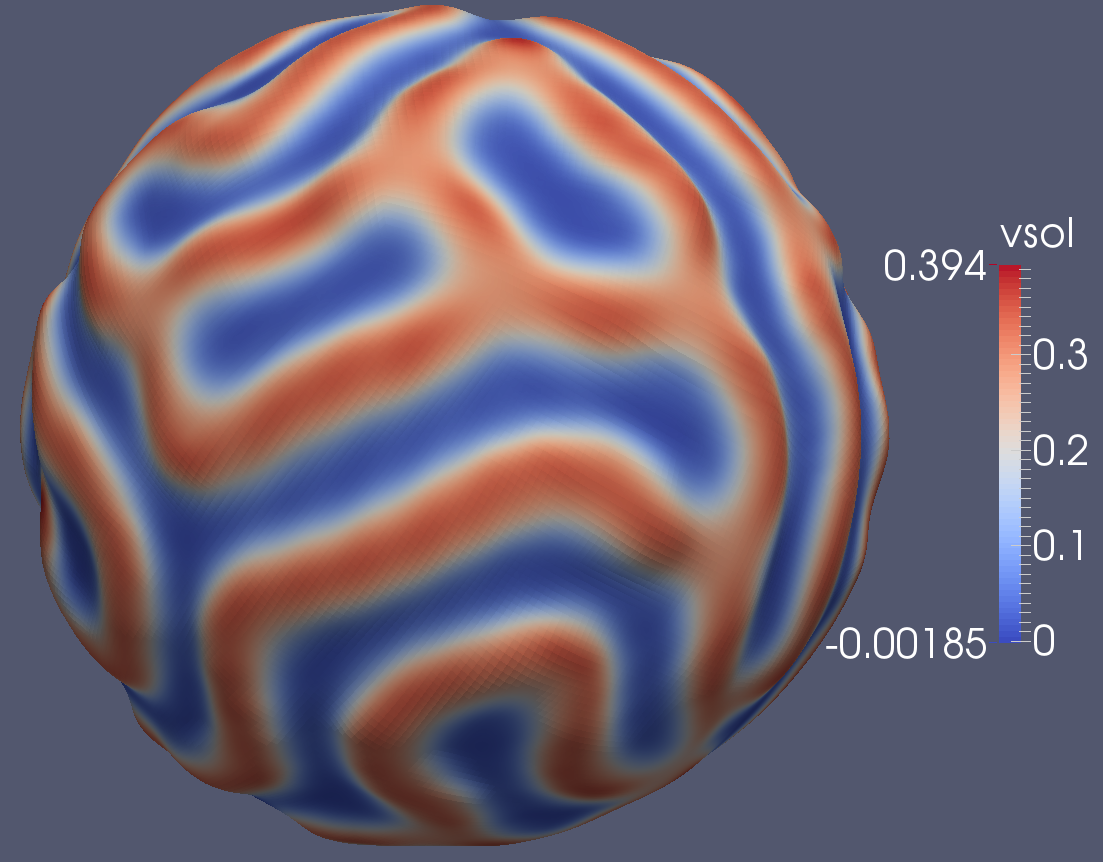}
    \caption{}
  \end{subfigure}
  \caption{Substrate $V$ and the geometry at $t^{400} = 1650$ (a) and $t^{800} = 6980$ (b).}
\label{fig:Sphere_1_004}
\end{figure}
\begin{figure}[H]
\centering
    \includegraphics[width=0.5 \textwidth]{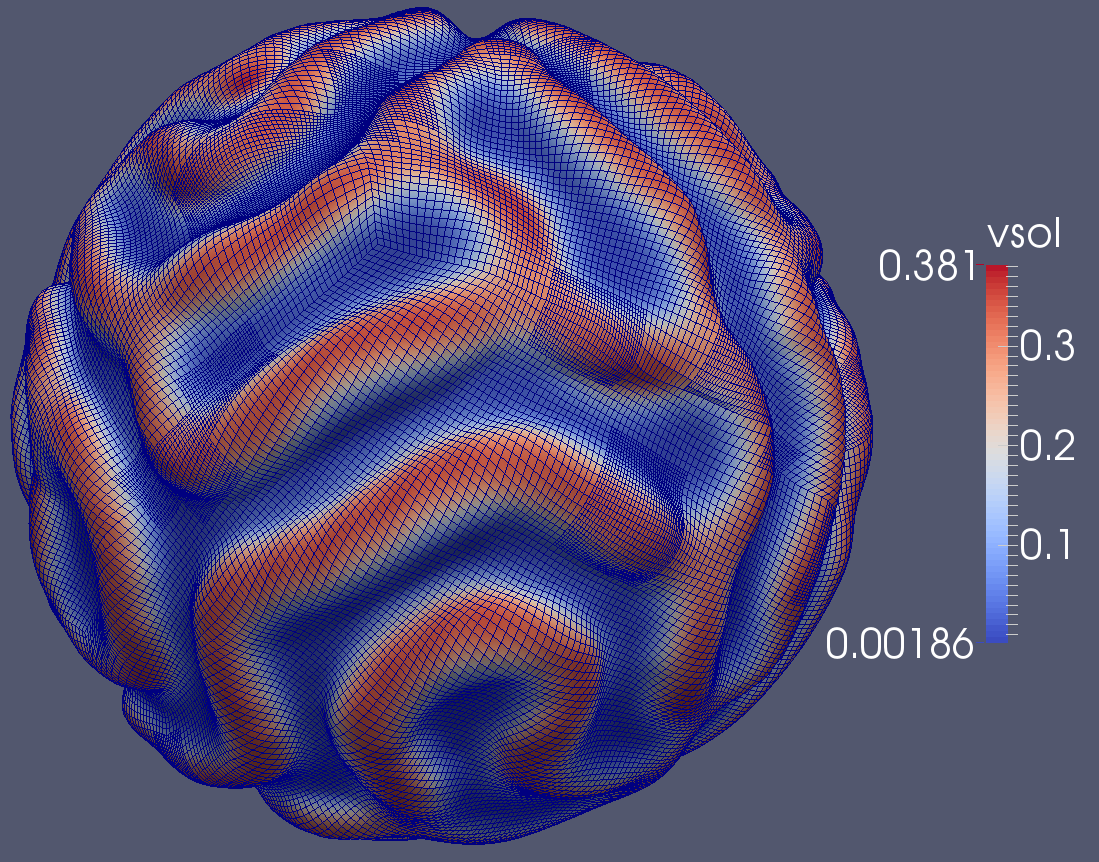}
  \caption{Substrate $V$ and the geometry at $t^{1200} = 1.3 \times 10^4$.}
\label{fig:Sphere_2_004}
\end{figure}
\begin{figure}[H]
\centering
  \begin{subfigure}[b]{0.45\textwidth}
    \centering
    \includegraphics[height=5cm]{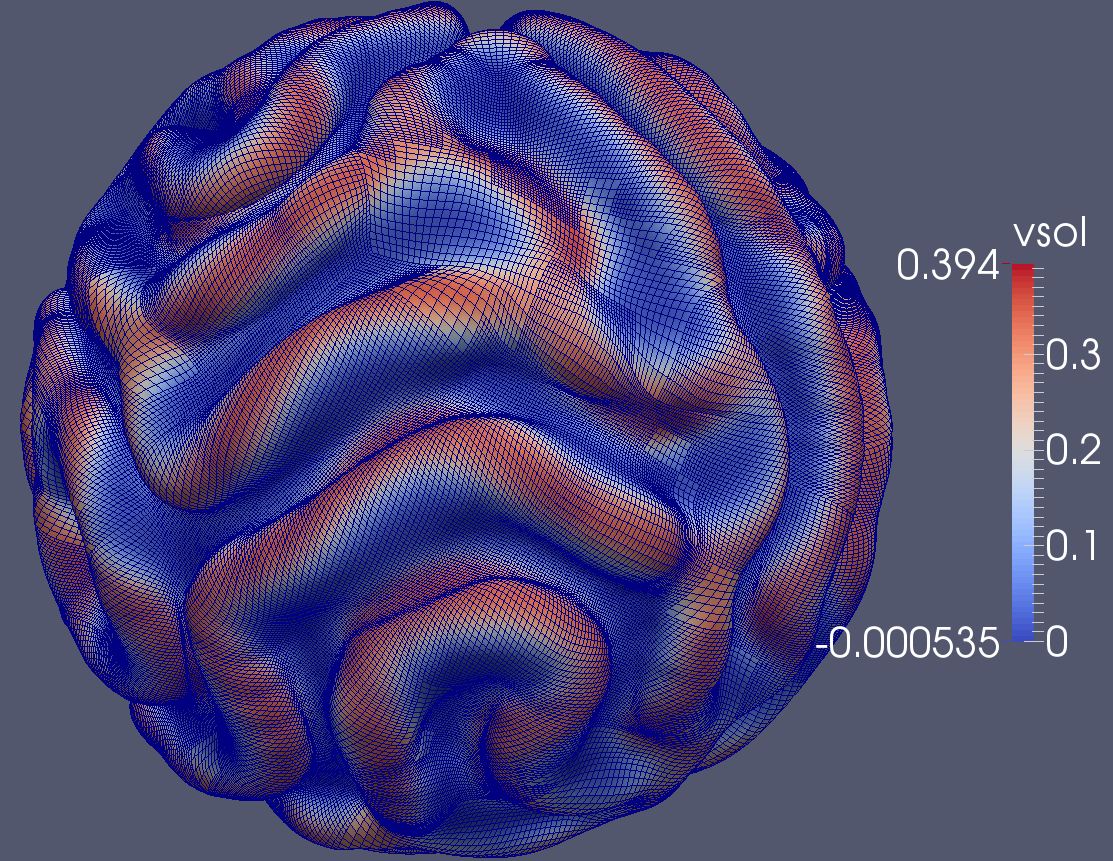}
    \caption{}
  \end{subfigure}
  \quad
  \begin{subfigure}[b]{0.45\textwidth}
    \centering
    \includegraphics[height=5cm]{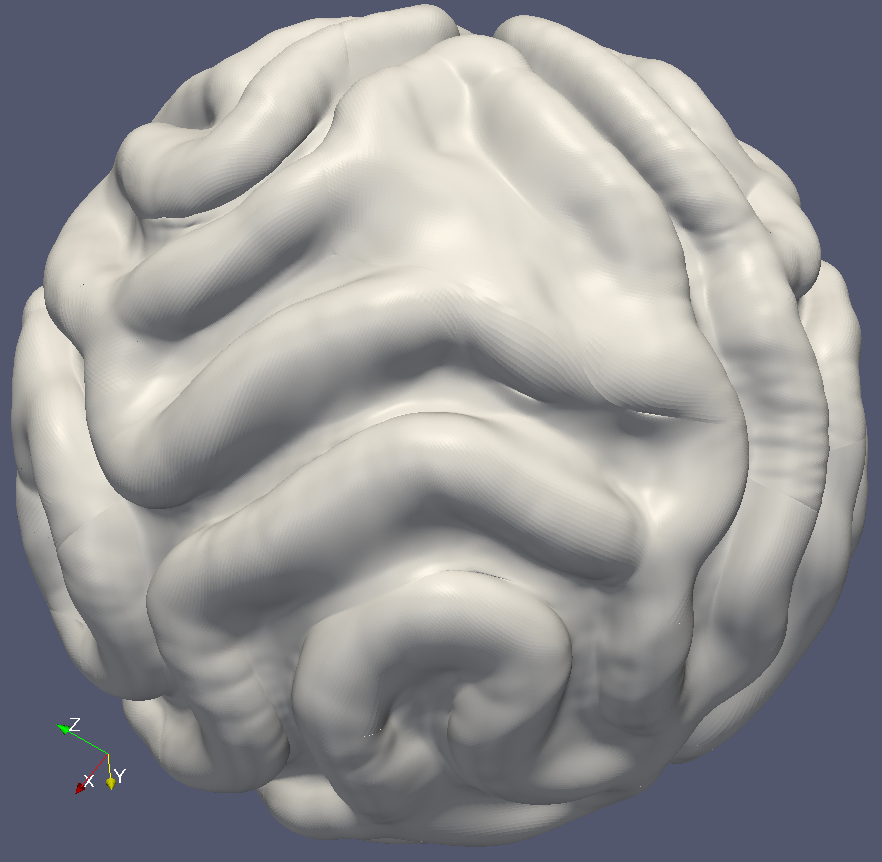}
    \caption{}
  \end{subfigure}
  \caption{Substrate $V$ and the geometry at $t^{1600} = 1.9 \times 10^4$ (a) and the final geometry for $F = 0.04$ at $t^{1650} \simeq 2 \times 10^{4}$ (b).}
\label{fig:Sphere_40_80_004}
\end{figure}

\noindent For a simulation with $(F, H, K, d_1, d_2) = (0.04, 0.06, 0.001, 0.2, 0.1)$, Figures \ref{fig:Sphere_1_004} to \ref{fig:Sphere_40_80_004} show the state in which concentration $V^k$ and the geometry $\mathcal{M}_k$ find themselves for various $t^k$. In Figure \ref{fig:Sphere_1_004} (a), we see that the four Gaussians have formed several narrow bands of nonzero concentration and in (b) we see strong pattern formation as well as the first signs of geometrical deformations. These deformations intensify in the course of time and then, the solver has performed the first local refinements in Figure \ref{fig:Sphere_2_004}, most likely due to curvature. At the later times, the solver has refined a large portion of the grid (see Figure \ref{fig:Sphere_40_80_004}). \\
The simulation has been terminated after $1650$ iterations. Slight unphysical undershoots in the concentrations are most likely due to truncation errors resulting from the discretization. For small time-steps $h_k$, the system matrix of the concentrations may lose its M-matrix property, due to a dominating presence of the mass matrix $[A]$ in equations (\ref{eq:System_u}) and (\ref{eq:System_v}). However, we have observed the undershoots (amounting to no more than $0.5$ percent of the characteristic magnitude of the concentrations) to be stable. We propose an improved numerical scheme, designed to suppress undershoots, in section \ref{subsect:undershoots}.

\subsection{Discussion}
\begin{figure}[H]
  \begin{subfigure}[b]{0.45\textwidth}
    \centering
    \includegraphics[height=6cm]{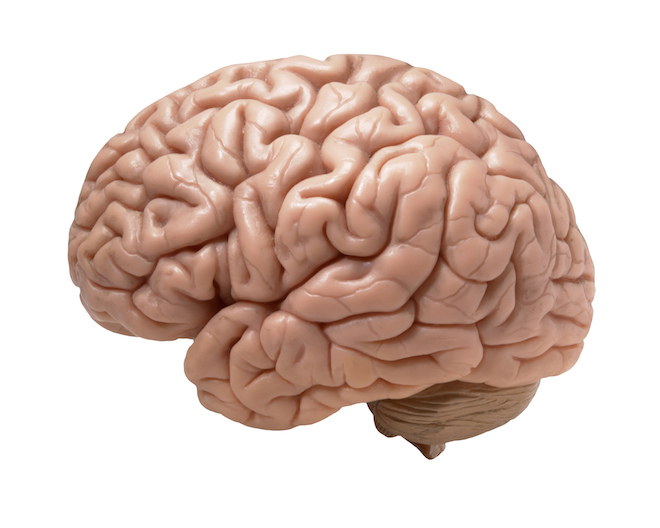}
    \caption{}
    \label{}
  \end{subfigure}
  \quad
  \begin{subfigure}[b]{0.45\textwidth}
    \centering
    \includegraphics[height=6cm]{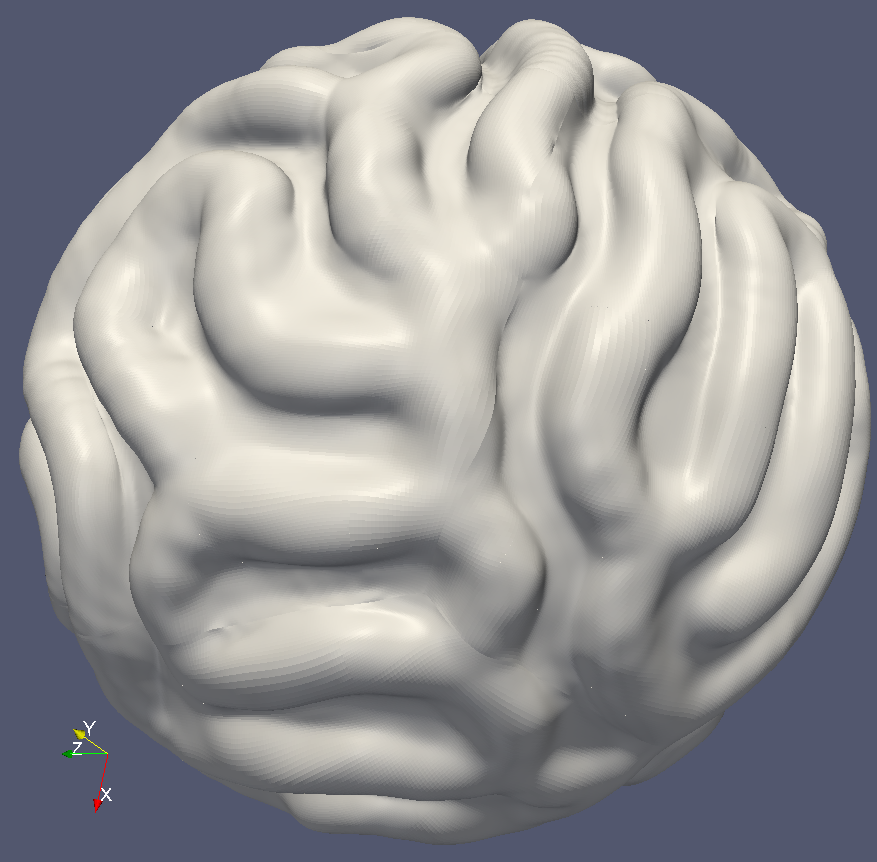}
    \caption{}
    \label{}
  \end{subfigure}
  \caption{Typical neural patterns in a healthy brain (a) and the results from the simulation with $F = 0.04$ (b).}
\label{fig:Results_Healthy_Brain}
\end{figure}

\begin{figure}[H]
\centering
  \begin{subfigure}[b]{0.45\textwidth}
    \centering
    \includegraphics[height=5cm]{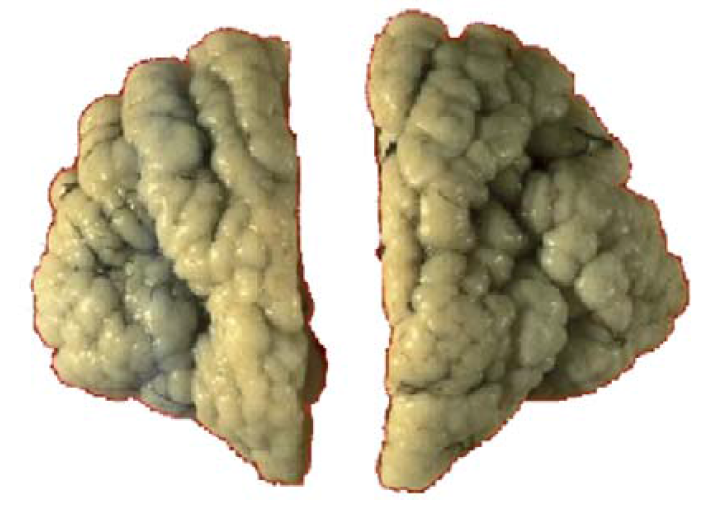}
    \caption{}
    \label{}
  \end{subfigure}
  \quad
  \begin{subfigure}[b]{0.45\textwidth}
    \centering
    \includegraphics[height=5cm]{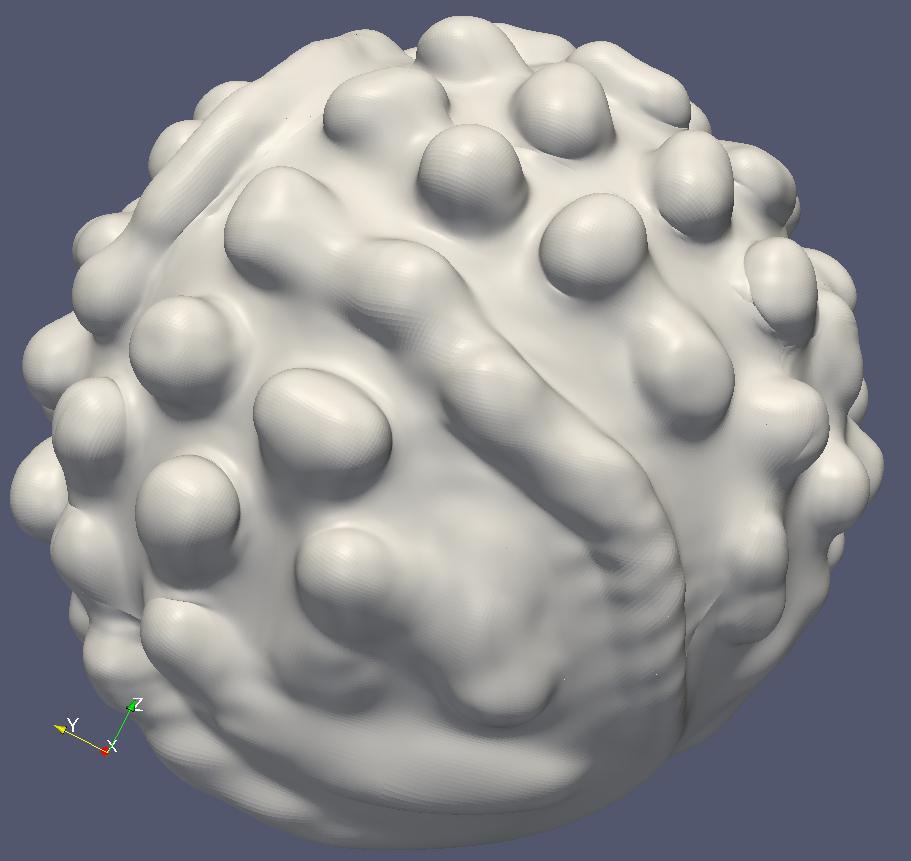}
    \caption{}
    \label{}
  \end{subfigure}
  \caption{Polymicrogyria (a) and the results from the simulation with $F = 0.0285$ (b).}
\label{fig:Results_Unhealthy_Brain}
\end{figure}

\begin{figure}[H]
\centering
  \begin{subfigure}[b]{0.45\textwidth}
    \centering
    \includegraphics[height=6cm]{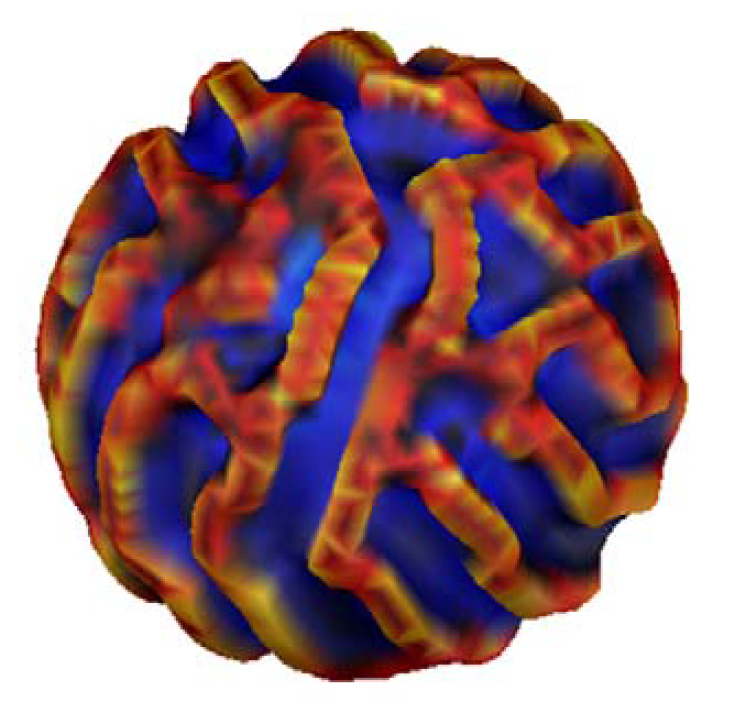}
    \caption{}
    \label{}
  \end{subfigure}
  \quad
  \begin{subfigure}[b]{0.45\textwidth}
    \centering
    \includegraphics[height=6cm]{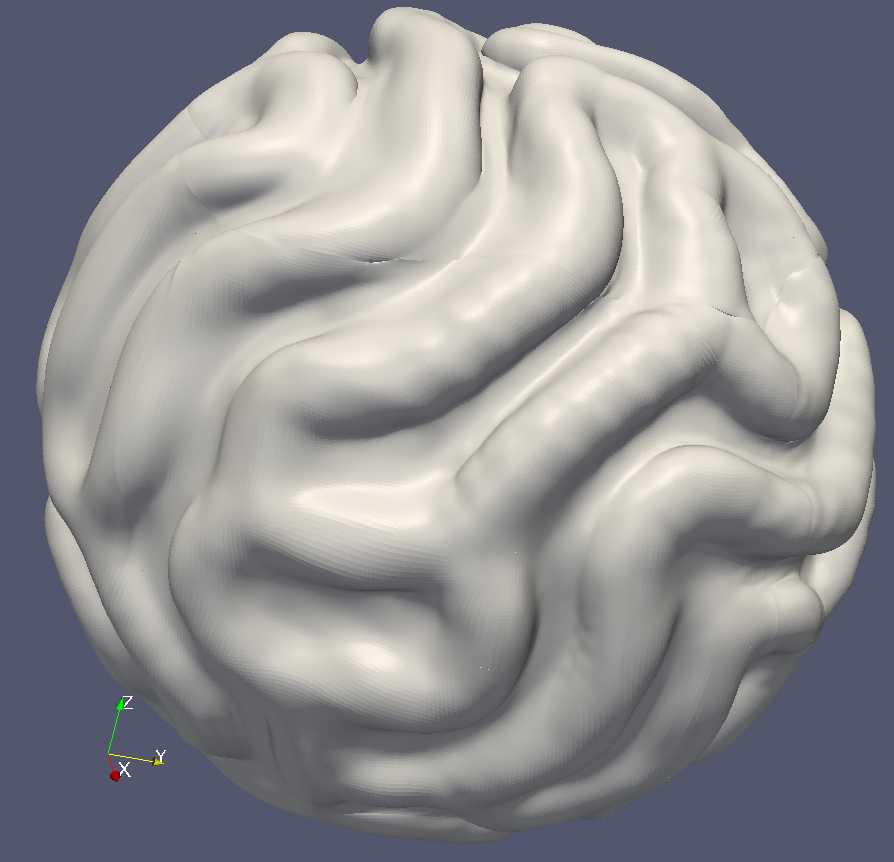}
    \caption{}
    \label{}
  \end{subfigure}
  \caption{The resulting geometry of a simulation with the classical FEM-approach from \cite{lefevre2010reaction} (a) and the geometry resulting from the IgA-scheme introduced in chapter \ref{chap:Isogeometric_Implementation} (b).}
\label{fig:Comparison_Results_Paper_IgA}
\end{figure}

The numerical scheme from Section \ref{chap:Isogeometric_Implementation} has been succesfully implemented on a six-patch topology and the results meet the expectations. We observe the characteristical pattern formation of the concentrations and the patterns manifest themselves in surface deformations. The deformations show a high degree of resemblance to typical brain patterns found in healthy adult individuals for $(F, H, K, d_1, d_2) = (0.04, 0.06, 0.001, 0.2, 0.1)$ (see Figure \ref{fig:Results_Healthy_Brain}).
\noindent Another implementation with the value of $F$ changed to $F = 0.0285$ shows mild resemblence with the neuropathology \emph{polymicrogyria} (see Figure \ref{fig:Results_Unhealthy_Brain}), which suggests that within the framework of this model, neuropathologies can be explained by deviations of the reaction rate $F$ due to genetical anomalies or other extrinsic influences.
\noindent With the above simulations, we were largely able to confirm the findings presented in \cite{lefevre2010reaction}. Note that we have used a growth factor $K = 0.001$ as opposed to $K = 0.0005$ in \cite{lefevre2010reaction}. Overall the results from the IgA-scheme, not surprisingly, exhibit improved smoothness when compared to their FEM-counterpart from \cite{lefevre2010reaction} (see Figure \ref{fig:Comparison_Results_Paper_IgA}) and the smoothness greatly contributes to the overall visual appeal of the results. It is noteworthy that the time-scale necessary to achieve similarly-sized folds as in \cite{lefevre2010reaction} is about a factor five larger even though the growth factor is doubled. A possible explanation is that we used a different initial condition.
\noindent There can be numerical artefacts due to the division of the spherical domain into patches. At the interfaces
between pairs of patches 'kinks' can arise due to reduced local smoothness of the IgA basis. The
'kink' may lead to oppositely directed normal vectors, which in turn may lead to unphysical various defects 
whenever the simulation is not terminated in time. These results can be improved by the application of
local smoothing of either the geometry or the normal vector in the vicinity of the interfaces between
adjacent patches. Smoothness across the patch boundaries can be improved by using the basis from \cite{Hinz2016Delft} (chapter $9$). The improved smoothness will most likely lead to greater visual appeal of the results but also to stronger folds as it may be possible to perform a larger amount of iterations without geometric clashing. Extending the model with a mechanism to avoid geometrical self-intersection by taking into account the mechanical stresses that occur during gyrification is difficult. However, it would most-likely lead to even more realistic outcomes. The effect of the patch boundaries on the concentrations has been studied in detail in \cite{Hinz2016Delft} (chapter Numerical Experiments).

\subsection{Suppressing Unphysical Undershoots}
\label{subsect:undershoots}
In Section \ref{chap:Results}, we noticed slight unphysical undershoots in the concentrations resulting from the fully discretized numerical scheme. In the following, we present a possible adjustment designed to prevent undershoots. In Section \ref{chap:Isogeometric_Implementation}, we presented a numerical scheme that is comprised of the following two main steps:
\begin{enumerate}
\item Compute the concentrations $u^{k+1}$ and $v^{k+1}$.
\item Update the mapping operator using $v^{k+1}$.
\end{enumerate}
The concentrations $u^k$ and $v^k$ are characterized by the corresponding vectors of weights $\mathbf{c}^k$ and $\mathbf{d}^k$, respectively. Let $\mathbf{x}_k = \left( \mathbf{c}^k, \mathbf{d}^k \right)$. We showed in Section \ref{chap:Results} that $\mathbf{x}_{k+1}$ is the result of solving a linear equation of the form
\begin{align}
\label{eq:main_equation_concentrations}
[Q] \mathbf{x}_{k+1} = \mathbf{f},
\end{align}
with $[Q] \in \mathbb{R}^{2N \times 2N}$ a positive-definite and block-diagonal system matrix. As such, equation (\ref{eq:main_equation_concentrations}) may be regarded as a quadratic optimization problem
\begin{align}
\label{eq:main_equation_concentrations_optimization}
\frac{1}{2} \mathbf{x}_{k+1}^T [Q] \mathbf{x}_{k+1}  - \mathbf{x}_{k+1}^T \mathbf{f} \longrightarrow \min \limits_{\mathbf{x}_{k+1}}.
\end{align}
A sufficient (however, for $p > 1$ not necessary) condition for the positivity of $u^{k}$ and $v^{k}$ is $\mathbf{x}_k \geq \mathbf{0}$. As such, we may augment the problem from (\ref{eq:main_equation_concentrations_optimization}) with this linear constraint, leading to a problem of the form
\begin{align}
\label{eq:main_equation_concentrations_optimization_augmented}
\min \limits_{\mathbf{x}_{k+1}} \quad & \frac{1}{2} \mathbf{x}_{k+1}^T [Q] \mathbf{x}_{k+1}  - \mathbf{x}_{k+1}^T \mathbf{f} \nonumber \\
& \text{s.t.} \quad \mathbf{x}_{k+1} \geq \mathbf{0}.
\end{align}
Problem (\ref{eq:main_equation_concentrations_optimization_augmented}) can be solved efficiently with standard convex optimization routines. Adding a constraint to the problem formulation is expected to increase the error per time-step in the norm induced by the matrix $[Q]$. However, we expect this additional error to be dwarfed by the truncation errors induced by the temporal and spatial discretizations. Heuristically, the vast majority of the entries in $\mathbf{x}_k$ are positive, even without constraints. The additional error introduced can be compensated for by local refinement of the basis.
\section{Conclusions}
\label{chap:Conclusions}
We have extended an existing numerical implementation of the Gray-Scott reaction-diffusion equations for surface deformation with the principles of \emph{Isogeometric Analysis}. This implementation holds in the most general way and has been applied to a multi-patch domain (sphere). Valid bases have been constructed by utilizing tensor-product B-splines in conjunction with the coupling of degrees of freedom. The numerical scheme has been enhanced with local refinement strategies and an adaptive time step selection. \\
We presented numerical results based on, among other parameters, two different choices for the values of the feeding rate $F$. For $F= 0.04$, we concluded that the results show strong resemblance with the characteristic patterns found on the surface of healthy adult human brains, adding credibility to the underlying reaction-diffusion model. For $F = 0.0285$, the results show mild resemblance to the neuropathology polymicrogyria, from which we concluded that, within the framework of this model, pathologies may be due to deviations of the underlying parameters. \\
The IgA-based results show a larger visual appeal compared to an equivalent FEM-implementation thanks to the smooth nature of the geometry description. It may be concluded that spline-based approaches are a viable alternative to classical FEM approaches whenever the quality of a model is assessed based on the quality of the numerical results, as in this case. \\
Regarding further research on this topic, basic space-time Galerkin schemes could be assessed.

  \end{document}